\documentstyle[amssymb,amsfonts]{amsart}

\newcommand{\K}{\mbox{{\bf K}}}  % disable Bbb
\def\be#1{\begin{equation} \label{#1}}
\def\bi{\begin{itemize}}
\def\bs{\begin{split}}
\def\es{\end{split}}
\def\ba{\begin{align}}
\def\bas{\begin{align*}}
\def\ea{\end{align}}
\def\eas{\end{align*}}
\def\R{{\hbox{\bf R}}}

\def\cube{{\hbox{\bf Q}}}

\def\sgn{{\hbox{sgn}}}
\def\dist{{\hbox{dist}}}
\def\R{{\hbox{\bf R}}}

\def\T{{\hbox{\bf T}}}

\def\eps{\varepsilon}

\newenvironment{proof}{\noindent {\bf Proof} }{\endprf\par}
\def \endprf{\hfill  {\vrule height6pt width6pt depth0pt}\medskip}
\def\emph#1{{\it #1}}
\def\textbf#1{{\bf #1}}

\def\Ga{\Gamma}

\def\Om{\Omega}

\def\pr{\partial}
\def\f12{\frac{1}{2}}
\newcommand{\nabb}{\mbox{$\nabla \mkern-13mu /$\,}}
\newcommand{\dabb}{\mbox{$D \mkern-13mu /$\,}}
\newcommand{\bea}{\begin{eqnarray}}
\newcommand{\eea}{\end{eqnarray}}
\newcommand{\beaa}{\begin{eqnarray*}}
\newcommand{\eeaa}{\end{eqnarray*}}

\theoremstyle{plain}
  \newtheorem{theorem}[subsection]{Theorem}
  
  \newtheorem{proposition}[subsection]{Proposition}
  \newtheorem{lemma}[subsection]{Lemma}
  \newtheorem{corollary}[subsection]{Corollary}

\theoremstyle{remark}

\theoremstyle{definition}
  \newtheorem{definition}[subsection]{Definition}

\include{psfig}

\begin{document}

\title[Bilinear estimates]{A physical space approach to wave equation bilinear estimates } 
\author{Sergiu Klainerman}
\address{Department of Mathematics, Princeton University, Princeton NJ 08544}
\email{ seri@@math.princeton.edu}

\author{Igor Rodnianski}
\address{Department of Mathematics, Princeton University, Princeton NJ 08544}
\email{ irod@@math.princeton.edu}

\author{Terence Tao}
\address{Department of Mathematics, UCLA, Los Angeles CA 90095-1555}
\email{ tao@@math.ucla.edu}

\dedicatory {Dedicated to the memory of Tom Wolff}

\vspace{-0.3in}
\begin{abstract}
Bilinear estimates for the wave equation in Minkowski space are normally proven 
using the Fourier transform and Plancherel's theorem.  However, such methods are
difficult to carry over to non-flat situations (such as wave equations with rough
metrics, or  connections with non-zero curvature).  In  this note we describe
an alternative  physical space approach which relies 
on vector fields, energy estimates as well as  tube localization, splitting into coarse
and fine scales, and induction on scales (in the spirit of Wolff \cite{wolff:cone},
\cite{wolff:smsub}). 
\end{abstract}

\maketitle

\section{Introduction}

In this paper $n \geq 2$ is an integer, and all implicit constants may depend on $n$.

We shall call a \emph{free wave} any (possibly vector-valued) solution\footnote{To avoid technicalities we shall assume that our waves are Schwartz functions in space.} $\phi$ to the free wave equation 
$$\Box \phi(x,t) = (-\frac{\partial^2}{\partial t^2} + \Delta) \phi(x,t) = 0$$ 
in $\R^{n+1}$.  We write $\phi[t] := (\phi(t), \phi_t(t))$ for the data of $\phi$ at $t$.  This space is equipped with the energy inner product
$$ \langle \phi[t], \psi[t] \rangle_e := \frac{1}{2} \int_{\R^n} \nabla 
\phi(x,t) \cdot \nabla \overline\psi(x,t) + \phi_t(x,t) \overline \psi_t(x,t)\
dx,$$
and we write $E(\phi[t]) := \langle \phi[t], \phi[t] \rangle_e$.
If $\phi$, $\psi$ are free waves, then these quantities are 
time-independent and we shall simply write them as $\langle \phi, \psi\rangle_e$ and
$E(\phi)$ respectively. 

If $\phi$ is a free wave and $\lambda$ is a dyadic number (i.e. an integer power of two), we say that $\phi$ \hbox{ has frequency $\lambda$} if the spatial Fourier transform $\hat \phi(\xi,t)$ is supported in the region\footnote{Here and in the sequel, we
 use $C$ to denote various constants depending only on $n$ and a parameter $\eps$ which
will appear later in the paper, and we use $A \lesssim B$ to denote the estimate $A \leq
CB$.  Finally, we use $A \sim B$ to denote the estimate $A \lesssim B \lesssim A$.}
$|\xi| \sim \lambda$; note that this property is independent\footnote{ The time
 independence property is invalid  in a curved spacetime.  There are however
 techniques, based on paradifferential calculus,
which allow one to circumvent this difficulty, see \cite{klainerman:quasil} 
and \cite{Kl-Rod}. } of
$t$.  

Let $\lambda$ be any dyadic number.  We use the symbol $P_\lambda$ to denote any Littlewood-Paley  projection operators of the form 
$$ \widehat{P_\lambda \phi}(\xi,t) := \eta(\xi/\lambda) \hat \phi(\xi,t)$$
where $\eta$ is any bump function adapted to the annulus $|\xi| \sim 1$; the exact choice of $\eta$ may vary from line to line.  Note that if $\phi$ is a free wave, then $P_\lambda \phi$ is a free wave with frequency $\lambda$.  Also, if $\phi$ has frequency $\lambda$, then we have the reproducing formula
\be{reproducing}
\nabla_x^b \phi = \lambda^{b-a} P_\lambda \nabla_x^a \phi
\end{equation}
for any $a, b \geq 0$ and some suitable $P_\lambda$ (possibly tensor-valued).
Here $\nabla_x^a$ refers to fractional spatial derivatives.
By computing the convolution kernel of $P_\lambda$ we thus obtain the pointwise estimates
\be{pointwise}
|\nabla_x^b \phi(x_0)| \lesssim \lambda^{b-a} \int |\nabla_x^a \phi(x_0 + y/\lambda)| (1 + |y|)^{-M}\ dy
\end{equation}
for any $M \geq 0$, with the implicit constant depending on $a$, $b$, $M$.  (Similar estimates obtain for spacetime derivatives $\nabla_{x,t}$ by using the equation $\Box \phi = 0$).
In particular we have the Sobolev inequality
\be{sobolev}
\| \phi \|_\infty \lesssim \lambda^{n/2-1} E(\phi)^{1/2}
\end{equation}
whenever $\phi$ is a free wave of frequency $\lambda$.  In practice, when we 
wish to estimate a bilinear expression $Q(\phi,\psi)$ in $L^2$, we shall put the lower
frequency wave in $L^\infty$ (using \eqref{sobolev} or more sophisticated variants), and
the higher frequency
 wave in $L^2$ (using energy estimates, possibly involving vector fields).

In the study of non-linear wave equations it has been realized 
that bilinear estimates on free waves play a fundamental role.  Two sample 
estimates are
given below; for a more thorough discussion of such estimates see
\cite{klainerman:nulllocal}, \cite{klainerman:foschi}, 
\cite{tao:hsdelta}, \cite{selberg}, \cite{ks}.

\begin{theorem}\label{null-1}    Let $\lambda$, $\mu$ be dyadic numbers, and let $\eps > 0$.  Then we have the estimate
\be{null-eq}
\begin{split}
\| Q(\phi, \psi) \|_{L^2(\cube_R)} \lesssim 
& R^\eps \min(\lambda,\mu)^{\frac{n-1}{2} + \eps}
E(\phi)^{1/2} E( \psi)^{1/2} \quad \hbox{ when }\,\, R \min(\lambda, \mu) \gtrsim 1\\
& R^{1/2} \min(\lambda,\mu)^{\frac{n}{2}}
E(\phi)^{1/2} E( \psi)^{1/2} \quad \hbox{ when }\,\, R \min(\lambda, \mu) \lesssim 1
\end{split}
\end{equation}
whenever $\cube_R$ is a spacetime cube of some side-length $R > 0$, 
$\phi$ and $\psi$ are free waves with frequency $\lambda$ and $\mu$ respectively, and $Q$ is one of the \emph{null forms}
\beaa
 Q_0(\phi,\psi) :&=& \phi_t \psi_t - \nabla \phi \cdot \nabla \psi,\\
Q_{\alpha \beta}(\phi,\psi) :&=& \partial_\alpha \phi \partial_\beta \psi - \partial_\beta \phi 
\cdot \partial_\alpha \psi.
\eeaa
\end{theorem}

\begin{theorem}\label{null-2}  Let $\mu, \lambda$ be dyadic numbers, and let
 $\phi, \psi$ be free waves of frequency $\lambda$.  Then we have
$$\| P_\mu( \phi \psi ) \|_{L^{q/2}_t L^{r/2}_x(\R^{n+1})} \lesssim 
(\frac{\mu}{\lambda})^{ n - \frac{4}{q} - \frac{2n}{r} }
\lambda^{n - 2 - \frac{2}{q} - \frac{2n}{r} }
E(\phi)^{1/2} E(\psi)^{1/2}$$
whenever $(q,r)$ are an \emph{admissible pair of Strichartz exponents} in the sense that
$$ 2 \leq q, r \leq \infty; \quad \frac{1}{q} + \frac{(n-1)/2}{r} \leq \frac{(n-1)/2}{2}; \quad (n,q,r) \neq (3,2,\infty).$$
\end{theorem}

Apart from the $\eps$ loss in Theorem \ref{null-1}, all the exponents
 are sharp \cite{klainerman:foschi}.  The $\eps$ can be removed; we will discuss this
issue later in the paper.

These two theorems apply only to frequency-localized free waves, but one can
 use Littlewood-Paley  decomposition and orthogonality in a standard manner to extend
these theorems to arbitrary free waves.  More precisely, we have

\begin{corollary}[``First generation'' null form estimate, see
\cite{klainerman:nulllocal}]\label{cor-1} If $Q$ is a null form and $\phi$, $\psi$ are
free waves, then
$$ \| Q(\phi, \psi) \|_{L^2_{t,x}([0,1] \times \R^n)} \lesssim
E(\langle D\rangle^{\frac{n-1}{2}+\eps}
\phi)^{1/2} E(\psi)^{1/2}
$$
for any $\eps > 0$, where $\langle D \rangle$ is the Fourier multiplier with symbol $(1 + |\xi|^2)^{1/2}$.
\end{corollary}

\begin{corollary}[Bilinear improvement to Strichartz' inequality]\label{cor-2}
Let $\phi$, $\psi$ be free waves, and let $(q,r)$ be an admissible pair of Strichartz exponents.  Then we have
$$
\| |D|^{-\sigma}( \phi \psi ) \|_{L^{q/2}_t L^{r/2}_x(\R^{n+1})} \lesssim 
E(|D|^{\frac{n-\sigma}{2} - \frac{1}{q} - \frac{n}{r}} \phi)^{1/2} E(|D|^{\frac{n-\sigma}{2} - \frac{1}{q} - \frac{n}{r}} \psi)^{1/2}$$
whenever $\sigma < n - \frac{4}{q} - \frac{2n}{r}$, and $|D|$ is the Fourier multiplier with symbol $|\xi|$.
\end{corollary}
The derivation of these corollaries from the Theorems is given in the Appendix.
We remark that Corollary \ref{cor-1} was
proven in \cite{klainerman:nulllocal} (with the sharp endpoint $\eps = 0$) by methods
heavily reliant on the Fourier transform.  Corollary \ref{cor-2} was similarly proven
using Fourier transform methods in \cite{kl-tat}, see also \cite{Kl:Ma2}.

The purpose of this note is to give different proofs of
 Theorems \ref{null-1}, \ref{null-2} and their corollaries
 which rely on physical space methods
instead of the Fourier transform.  As such, we believe these
 arguments will be useful for
such situations as quasilinear wave equations with rough data, in which the Fourier
transform is difficult to use directly.  The techniques used here extend to many other
types of bilinear estimates, but we shall restrict ourselves to the above two types of
estimates for simplicity.

We briefly outline the ideas used in our arguments:

\begin{itemize}

\item \textbf{Splitting physical space into coarse scales and fine scales.}  If a 
large cube $\cube_R$ is partitioned into disjoint\footnote{We always ignore sets of
measure zero when we say that two sets are disjoint.} sub-cubes $\cube_r$, we can split 
\be{p-split}
\| \phi \|_{L^p(\cube_R)} = 
(\sum_{\cube_r: \cube_r \subset \cube_R} \| \phi \|_{L^p(\cube_r)}^p)^{1/p}.
\end{equation}
The inner $L^p(\cube_r)$ norm represents the fine scales, while the
 outer $l^p$ norm represents the coarse scales.  The strategy is then to choose $r$
carefully and handle the coarse and fine scales in different ways.  This should be
compared with Fourier-based methods, which often employ a similar decomposition but in
the frequency variable (\cite{kl-tat}, \cite{tao:hsdelta}, etc.).

\item \textbf{Decomposition into wave packets.}  In a spacetime 
cube $\cube_R$, one can split a free wave $\phi$ into almost orthogonal ``wave packets''
$\phi_T$, each of which are also free waves and are essentially localized in a tube $T$,
which is a suitable neighbourhood of a light ray.  The traditional decomposition
into wave packets, which has its origin in  the so called second
 dyadic decomposition\footnote{These were developed further by Bourgain
\cite{borg:kakeya} and later authors (\cite{tvv:bilinear}, \cite{tv:cone1},
\cite{wolff:cone}, \cite{tao:cone}, \cite{wolff:smsub}, etc.)} pioneered by 
 Fefferman \cite{Feff}
 and C\'ordoba \cite{Cord}
in connection to the $L^p$ multiplier problem, relies heavily on the  Fourier
transform. 

In our paper we suggest a different construction of wave packets
based on a ``dynamical'' physical space decomposition 
which   avoids Fourier space decompositions 
and   the related  method of ``stationary phase''.
This is achieved by localizing our waves at two distinct times.
One can use the vectorfield method, instead of stationary phase, to
show that the resulting waves concentrate along tubes.

\item \textbf{Distinguishing between transverse and parallel interactions.} 
 Once $\phi$ and $\psi$ are decomposed into wave packets, one can then separate into two
cases depending on the angle between the two packets.  When the angle is small,
 we say we
have a \emph{parallel interaction}, and we expect the cancellation 
in the null structure
to play a role. 
 In the
large angle case we have a \emph{transverse interaction}, which we expect to be well
localized in physical space. In the original approach to the proof of
 bilinear estimates the two cases\footnote{or more precisely
their Fourier space counterparts.} were both treated with
the help of the  spacetime  Fourier transform. 
In some cases one can use Lorentz transforms to 
reduce the
parallel case to the transverse case: see \cite{tvv:bilinear}, \cite{tv:cone2},
\cite{wolff:cone}, \cite{wolff:smsub} for further discussion.

The new twist in our paper is the treatment of parallel
interactions, for which the vectorfield method turns out  to be very well
suited. Together with the dynamic tube localization and the treatment
of transversal interactions, based on induction on scales initiated in  
\cite{borg:kakeya}, \cite{tvv:bilinear},
\cite{tv:cone2},
\cite{wolff:cone}, \cite{wolff:smsub}, this allows us to give an 
almost\footnote{With the exception of the dyadic
 decomposition.} completely geometric
treatment of the first bilinear estimate, see Theorem \ref{null-1}

\item \textbf{Induction on scales (or bootstrap arguments).}  
To prove an estimate at a large spatial scale $R$, assume inductively that an estimate
has already been obtained at smaller scales $r$.  This idea has been around for some time
(see e.g. \cite{borg:kakeya}, \cite{tvv:bilinear}, \cite{tv:cone1}), but only in the
recent work of Wolff \cite{wolff:cone}, \cite{wolff:smsub} was it realized that such
techniques can give very sharp estimates.

\item \textbf{Improved decay estimates.} 
The original proof of Theorem \ref{null-2} was based
on an improved $L^1-L^\infty$  decay estimate which   takes 
advantage of the size of the
Fourier support of the data in a small rectangle in 
Fourier space.

 Here we take a different approach, first developed 
by \cite{tao:lowreg}, based on a similar phenomenon
in physical space. This  hinges on the observation
that  ``bump waves'', i.e.
solution of the free wave equation with initial data localized at frequency one
and decreasing rapidly away from the unit ball centered at the origin, are 
essentially supported on a set of codimension one. 
 In flat space this can be easily seen
 from the explicit representation of  bump waves as 
Fourier integrals. Due to this fact one
can
 improve the
$L^\infty$ decay of bump waves  by averaging over a
relatively coarse scale in physical space. 
Here we show that this improved decay property can be
derived by vectorfield techniques in the spirit of \cite{klainerman:quasil}.
Thus, by avoiding  the explicit representation of solutions, we give
an essentially geometric derivation of Theorem \ref{null-2}. This
opens
the way, we believe, for significant generalizations to 
quasilinear equations\footnote{  The recent result of \cite{Bahouri:Chemin}
is based on a generalization of the bilinear estimates in Theorem \ref{null-2}
to wave equations with rough coefficients. The construction is based
on Fourier Integral operators. }.

 By frequency-localized fundamental solution we mean a free wave of frequency 1 whose
initial data is rapidly decreasing away from the origin.  It is well known that the
$L^\infty$ decay of such a solution in time can be combined with energy estimates to
yield Strichartz estimates.  However, there are other properties of the solution which
are also useful.  For instance, the fundamental solution is essentially supported on a
set of codimension 1, which means that one can improve the $L^\infty$ decay of the
(absolute value of) the fundamental solution by averaging over a relatively coarse scale
in physical space.  This fact (first used in \cite{tao:lowreg}) can be used to control
the low-frequency portion of a bilinear or nonlinear expression.  We remark that vector
fields methods are well suited for obtaining this type of control on the fundamental
solution, even in rough metrics (see e.g. \cite{klainerman:quasil}).

\end{itemize}

In the proof of Theorem \ref{null-1} we shall use the first four ideas, while in the proof of Theorem \ref{null-2} we shall use the first and fifth ideas.  To realize these ideas we shall use energy estimates and
commutation with special vectorfields in the spirit
 of \cite{klainerman:quasil}, see also \cite{klainerman1}, \cite{klainerman2}.  It is possible to use other methods (e.g. Fourier integral parametrices) as a substitute for vector field methods, but we believe the vector field approach is better adapted for more non-linear situations (e.g. quasilinear wave equations) in which the Fourier transform is less amenable.
This paper is organized as follows.  We prove Theorem \ref{null-1} in Sections \ref{tubes-sec} to \ref{second-proof}.  In Section \ref{tubes-sec} we set up
the basic machinery of tubes or wave packets, and in Section \ref{induction-sec}
we combine this with an induction on scales argument to reduce to the
parallel interaction case.  After some vector fields preliminaries in
Section \ref{vector-sec} we handle the parallel interaction case in 
Section \ref{parallel-sec}.  Finally, we give two proofs of the key tube
localization lemma (one a hybrid of spatial and Fourier localization methods,
and the other based entirely on spatial localization) in Sections \ref{first-proof} and \ref{second-proof}.  After some remarks in Section \ref{remarks-sec} we then prove Theorem \ref{null-2} in Section \ref{null-sec}.

The first author is supported by NSF grant DMS-0070696. The  second author is supported 
by NSF grant DMS-0107791.
The third author is a Clay Prize Fellow and is supported by a grant from the Packard Foundation.  This work was conducted at UCLA, Princeton, Oberwolfach, and ETH Zurich; the first and third authors are especially grateful for Michael Struwe's hospitality during their stay at ETH.

\section{Tube localization}\label{tubes-sec}

A fundamental tool in our proof of Theorem \ref{null-1} is the ability to decompose a free wave into ``wave packets'' which are concentrated in tubes along null rays.  To make this precise we need some notation.  We shall be working in the standard spacetime cube
\be{cube0-def}
\cube^{(0)}_R := \{ (x,t): 0 \leq x_i \leq R \hbox{ for all } i=1,\ldots, n; \quad R \leq t \leq 2R \}
\end{equation}
for some $R \gtrsim 1$.  Note that this cube is some distance away from the initial spatial hypersurface $t = 0$; this will be very convenient when we apply vector field techniques.

We shall also need a small constant $0 < \eps \ll 1$; all our implicit constants may depend on $\eps$.

\begin{definition}\label{tube-def}  Let $1 \lesssim r \lesssim R$ be scales, let $\omega \in S^{n-1}$ be a direction, and let $x_0 \in \R^n$ be a point in space.  We call the set
$$ T = T(R,r,\omega,x_0) := \{ (x,t): |x - t \omega - x_0| \leq r; 0 \leq t \leq 3R \}$$
the \emph{tube} with \emph{length} $R$, \emph{width} $r$, \emph{velocity} $\omega$, and \emph{initial position} $x_0$. 
\end{definition}

We now define the notion of a wave packet localized to a tube.  We shall
 normalize our definition in the situation when the original waves 
$\phi$, $\psi$ are normalized to
 have energy 1.  

\begin{definition}\label{packet-def}  Let $T$ be a tube and $\phi_T$ be
 a free wave with frequency $\mu \gtrsim 1$ and energy $E(\phi_T) \lesssim 1$. 
 We say that $\phi_T$ is a \emph{wave packet localized to $T$} if one has the
 estimates
\be{local}
| \nabla_{x,t} \phi_T(x,t) | \lesssim (R + \dist((x,t), T))^{-100n} \mu^{n/2}
\end{equation}
whenever $R/2 \leq t \leq 5R/2$ and $(x,t) \not \in T$, as well as the $L^2$ variant
\be{local-2}
(\int_{x: |x| \lesssim R, (x,t) \not \in T}
| \nabla_{x,t} \phi_T  |^2\ dx)^{1/2} \lesssim (R + \dist((0,t), T))^{-100n}
\end{equation}
for all $R/2 \leq t \leq 5R/2$.
\end{definition}

The estimate \eqref{local} essentially follows from \eqref{local-2} by \eqref{pointwise}, perhaps at the cost of replacing $T$ by a tube of twice the radius.  (The error term coming from inside the tube can be controlled by \eqref{sobolev}).  In practice we shall use \eqref{local} for low-frequency waves and \eqref{local-2} for high-frequency waves.

The main tube decomposition lemma can now be stated as follows.

\begin{lemma}[Tube decomposition lemma]\label{tube-decomp}  Let $\phi$ be a free wave with energy $E(\phi) = 1$ and frequency $\mu \gtrsim 1$, and let $R \gtrsim 1$ and $R^{1/2+\eps} \lesssim r \lesssim R$.  Then there exists a collection $\T$ of tubes of length $R$ and width $r$, 
and to each tube $T$ in $\T$ there exists a wave packet $\phi_T$ localized to $T$ with frequency $\mu$, such that we have the decomposition
\be{decomposition}
\phi = \sum_{T \in \T} \phi_T + \phi_{error}
\end{equation}
where the error term is a free wave which satisfies $E(\phi_{error}) \ll E(\phi)$.  Also, we have
the Bessel inequality
\be{bessel}
\sum_{T \in \T} E(\phi_T) \lesssim 
\end{equation}
and the almost orthogonality relationship
\be{ortho}
E(\sum_{T \in \T'} \phi_T) \lesssim R^{C\eps} \sum_{T \in \T'} E(\phi_T) + R^{-100n}
\end{equation}
for any subset $\T'$ of $\T$.  Finally, to control error terms we require the technical property that given any cube $\cube_R$, the number of tubes in $\T$ which intersect $\cube_R$ is $O(R^{10n})$.
\end{lemma}

At a heuristic level one should imagine that $\phi_{error} = 0$, and that the $\phi_T$ are an orthogonal decomposition of $\phi$.

This type of lemma has been used by many authors.  For the closely related restriction problem for the Fourier transform, these wave packets first appear in the work of Fefferman, Knapp, and C\'ordoba, and were developed substantially further by Bourgain \cite{borg:kakeya}, \cite{borg:cone}, Wolff \cite{wolff:cone}, \cite{wolff:smsub}, and others (e.g. \cite{tao:cone}, \cite{tataru:bilinear}).  In these papers the tube decomposition is achieved by a decomposition in both physical and frequency space.  We shall give two proofs of this lemma in Sections \ref{first-proof}, \ref{second-proof}.  The first proof is similar to the previous arguments but only uses the Fourier transform at the level of the initial data rather than in spacetime, and uses commutation with vector fields to control the localization.  The second proof is based purely on physical space methods, and achieves localization to a tube by localizing the wave to a disk at two different times (see  somewhat 
related  ``slices'' approach to the Kakeya problem in e.g. \cite{borg:high-dim}).

In the high frequency case $\mu \gg 1$ the condition $r \gtrsim R^{1/2+\eps}$ can be relaxed to $r \gtrsim R^{1/2+\eps} \mu^{-1/2+\eps}$, but we are unable to exploit this for our bilinear estimates because our argument needs $\phi$ and $\psi$ to be decomposed into tubes of equal width.

Assume Lemma \ref{tube-decomp} for the moment.  In the next section we see how this lemma can be combined with an induction on scales argument to reduce the bilinear estimate \eqref{null-eq} to the ``parallel interaction case''.  In the section after that we show how commutation with vector fields can be used to treat these parallel interactions, thus completing the proof of \eqref{null-eq}.

\section{Induction on scales}\label{induction-sec}

We now begin the proof of Theorem \ref{null-1}.  The idea is to use Lemma \ref{tube-decomp} and an induction on scales (or bootstrap) argument to reduce \eqref{null-eq} to the parallel interaction case.

We first give some preliminary reductions.  Fix $\lambda$ and $\mu$; without loss of generality we may assume that $\lambda \leq \mu$.  By rescaling we may take $\lambda = 1$.

The case $R \lesssim 1$ is easy to dispose of by H\"older and the Sobolev inequality \eqref{sobolev}:
\be{stupid}
\| Q(\phi, \psi) \|_{L^2(\cube_R)} \lesssim R \| \nabla_{x,t} \phi \|_{L^\infty_t L^\infty_x} \| \nabla_{x,t} \psi \|_{L^\infty_t L^2_x} 
\lesssim R^{\frac 12} E(\phi)^{1/2} E(\psi)^{1/2}.
\end{equation}
Thus we may take $R \gg 1$.

We now set up the bootstrap argument.
For any $\alpha > 0$, let $B(\alpha)$ denote the statement that the bilinear estimate 
\begin{equation}
\| Q(\phi, \psi) \|_{L^2(\cube_R)} \lesssim 
R^\alpha E(\phi)^{1/2} E( \psi)^{1/2}
\tag{$B_\alpha$}
\end{equation}
holds for all free waves $\phi$, $\psi$, all $R \gg 1$ with frequency 1 and $\mu$ respectively, and all spacetime cubes $\cube_R$; the implicit constant can depend on $\alpha$ and the implicit constants used to define such notions as ``$\psi$ has frequency $\mu$''.  Our objective is to prove $B(\alpha)$ holds for all $\alpha > 0$.

From \eqref{stupid} we have the estimate $B(\frac 12)$.  The claim will now follow from bootstrapping the following Proposition.

\begin{proposition}\label{bootstrap}  If $B(\alpha)$ holds for some $\alpha > 0$, then $B((1-c)\alpha + C\eps)$ holds for all $\eps > 0$, where $c>0$ is an absolute constant depending only on $n$.
\end{proposition}

Indeed, from Proposition \ref{bootstrap} it follows that the infimum of all $\alpha > 0$ for which $B(\alpha)$ holds must be zero.

It remains to prove the Proposition.  Fix $\alpha > 0$, $\eps > 0$, and assume that $B(\alpha)$ already holds.  Fix free waves $\phi$, $\psi$ and a scale $R \gg 1$, and a cube $\cube_R$.  We may assume that we have the normalization $E(\phi) = E(\psi) = 1$.  We shall take advantage of spacetime translation invariance to set $\cube_R = \cube^{(0)}_R$, where $\cube^{(0)}_R$ was defined in \eqref{cube0-def}; note in particular that the cube $\cube_R$ is a distance $R$ away from the initial spatial hypersurface $t = 0$.  Our objective is to prove that
\be{objective}
\| Q(\phi, \psi) \|_{L^2(\cube^{(0)}_R)} \lesssim 
R^{(1-c)\alpha + C\eps}.
\end{equation}
If one applies the hypothesis $B(\alpha)$ directly to this cube $\cube^{(0)}_R$ one only gets a bound of $R^\alpha$.  To improve upon this we shall only apply $B(\alpha)$ at smaller scales than $R$, and use transversality to sum efficiently (in the spirit of \cite{wolff:cone}).

We let $R^{1/2 + \eps} \ll r \ll R$ be an intermediate dyadic scale between $\sqrt{R}$ and $R$ to be chosen later; $r$ will eventually be fairly close to $\sqrt{R}$. 
We now use Lemma \ref{tube-decomp} to write
$$ \phi = \sum_{T \in \T} \phi_T + \phi_{error}; \quad \psi = \sum_{T' \in \T'} \psi_{T'} + \psi_{error}$$
where $\T$, $\T'$ are collections of tubes with length $R$ and width $r$, 
and for each $T \in \T$, $\phi_T$ is a wave packet localized to $T$ obeying
\eqref{bessel}, \eqref{ortho}.  Similarly for the $\psi_{T'}$.

The $\phi_{error}$, $\psi_{error}$ are small in energy and their contribution can be iterated away by the usual argument\footnote{Specifically, one assumes that $B(\alpha)$ holds with an explicit constant $A$, possibly depending on $R$, and use the small energy of $\phi_{error}$ and $\psi_{error}$ to bound the contribution of these terms by at most $A/10$.  At the end one can absorb the $A/10$ into the left hand side.}.  We can therefore split the remaining portions of the null form $Q(\phi,\psi)$ into a sum over interactions between wave packets:
\be{packet}
\sum_{T \in \T} \sum_{T' \in \T'} Q(\phi_T, \psi_{T'}).
\end{equation}

By \eqref{local} and H\"older one may easily 
dispose of those tubes $T$ for which $\dist(T, \cube^{(0)}_R) \gtrsim R$. 
 Similarly for $T'$.  Thus we shall assume in the sequel that $\dist(T, \cube^{(0)}_R), \dist(T', \cube^{(0)}_R) \lesssim R$ for all $T \in \T$ and $T' \in \T'$.  In particular the cardinality of $\T$ and $\T'$ is only $O(R^{10n})$.

We shall need an intermediate separation scale $r \ll \rho_0 \ll R$ to be chosen later; eventually $\rho_0$ will be close to $r$.  We will assume that $\rho_0$ is a dyadic number.

Informally, the strategy will be as follows.  We wish to divide the above sum into \emph{parallel interactions} $|x_0(T) - x_0(T')| \leq \rho_0$ and \emph{transverse interactions} $|x_0(T) - x_0(T')| \geq \rho_0$, where $x_0(T)$ denotes the initial position of $\,T$.  The parallel interactions should be small thanks to the cancellation in the null form, and we will use vector field methods to exploit this.  For the transverse interactions, we do not exploit the null structure, but instead we use transversality to ensure that $T \cap T'$ is localized to a cube of length $\ll R$, so that when we apply the hypothesis $B(\alpha)$ we shall gain in the $\alpha$ exponent.  Of course we still have to sum over all tubes $T$, $T'$, but it turns out that we have enough orthogonality in \eqref{bessel}, \eqref{ortho} to achieve this.

We now turn to the details.  For technical reasons we have to take a little more care in dividing the wave packet interactions into parallel and transverse interactions in order for the transverse interactions to decouple nicely (cf. \cite{tvv:bilinear}).  Specifically, we shall use a bilinear Whitney decomposition of the initial data.

The initial positions $x_0$ of the tubes in $\T$, $\T'$ can all be contained in a spatial cube $K$ of dyadic side-length $CR$.  Let $\K$ denote the set of all dyadic sub-cubes $\kappa$ of $K$ whose side-length $l(\kappa)$ is at least $\rho_0/2$. Without loss of generality we may assume that the initial positions $x_0$ of tubes in $\T$, $\T'$ never lie on the boundary $\partial \kappa$ of any cube $\kappa \in \K$.

If $\kappa, \kappa' \in \K$ are two such sub-cubes of the same side-length, we say that $\kappa$ and $\kappa'$ are \emph{close} if $\kappa$ is not adjacent to $\kappa'$, but the parent of $\kappa$ is adjacent to the parent of $\kappa'$.  We write $\kappa \sim \kappa'$ to denote the statement that $\kappa$ and $\kappa'$ are close.  Note that if $\kappa \sim \kappa'$, then $\dist(\kappa,\kappa')$ is comparable to the side-length of $\kappa$ or $\kappa'$.  

For each $\kappa \in \K$, define the sets 
$$ \T(\kappa) := \{ T \in \T: x_0(T) \in \kappa \}
\hbox{ and } \T'(\kappa) := \{ T' \in \T': x_0(T') \in \kappa \};$$
these represent the tubes emanating from the cube $\kappa$.  A simple stopping time argument reveals that for any $T \in \T$ and $T' \in \T'$, exactly one of the following statements is true:

\begin{itemize}
\item (Transverse interaction) There exists a close pair of cubes $\kappa \sim \kappa'$ such that $T \in \T(\kappa)$ and $T' \in \T'(\kappa')$.

\item (Parallel interaction) There exists cubes $\kappa$, $\kappa'$ of side-length $l(\kappa) = \rho_0$ which are either equal or adjacent such that $T \in \T(\kappa)$ and $T' \in \T'(\kappa')$.  (We shall write $\kappa \approx \kappa'$ in this case.
\end{itemize}

We can therefore split the null form \eqref{packet} into the \emph{transverse interaction component}
\be{transverse}
\sum_{\kappa, \kappa' \in \K: \kappa \sim \kappa'} Q( \phi_\kappa, \psi_{\kappa'} )
\end{equation}
and the \emph{parallel interaction component}
\be{parallel}
\sum_{\kappa,\kappa' \in \K: \kappa \approx \kappa', l(\kappa) = l(\kappa') = \rho_0} Q( \phi_\kappa, \psi_{\kappa'} )
\end{equation}
where
$$ \phi_\kappa := \sum_{T \in \T(\kappa)} \phi_T \hbox{ and }
\psi_\kappa := \sum_{T' \in \T'(\kappa)} \psi_{T'}.$$
Informally, $\phi_\kappa$ represents the portion of $\phi$ emanating from $\kappa$, and similarly for $\psi_\kappa$.

Consider first the contribution of the transverse component \eqref{transverse}.  Observe that there are only $O(\log R)$ possible values for the side-length $l(\kappa)$.  Since we are allowing ourselves to lose an $R^\eps$ in the final estimate \eqref{objective}, it suffices by the triangle inequality to show that
$$
\| \sum_{\kappa, \kappa' \in \K: \kappa \sim \kappa'; l(\kappa) = l(\kappa') = \rho} Q( \phi_\kappa, \psi_{\kappa'} )
\|_{L^2(\cube^{(0)}_R)} \lesssim R^{\alpha-\sigma+C\eps} 
$$
for some $\sigma>0$ and for all $\rho_0 \leq \rho \leq 1$.  

Fix $\rho$.  We use the triangle inequality to pull the sum in $\kappa, \kappa'$ outside the norm.  Observe that for each $\kappa$ there are only $O(1)$ cubes $\kappa'$ which are close to $\kappa$, and vice versa.  By Cauchy-Schwarz and Bessel's inequality \eqref{bessel} it thus suffices to show that
\be{transverse-bound}
\| Q( \phi_\kappa, \psi_{\kappa'} )
\|_{L^2(\cube^{(0)}_R)} \lesssim R^{\alpha-\sigma+C\eps} 
(\sum_{T \in \T(\kappa)} E(\phi_T))^{1/2} 
(\sum_{T' \in \T'(\kappa')} E(\psi_{T'}))^{1/2} + R^{-50n}
\end{equation}
for all close pairs of cubes $\kappa$, $\kappa'$ of side-length $\rho$.

Fix $\kappa$, $\kappa'$.
Observe that if two tubes $T$, $T'$ in the above expression intersect, the angle between $T$ and $T'$ is comparable\footnote{This crucial fact, that separation in space implies separation in angle, relies on the choice of $\cube^{(0)}_R$ to be away from the initial spatial hypersurface $t=0$.} to $\rho/R$.  Since $T, T'$ have width $r$ and $\rho \geq \rho_0 \gg r$, we thus see by elementary geometry that $T \cap T'$ is contained in a cube of side-length $O(Rr/\rho) \ll R$.  To exploit this we shall subdivide the cube $\cube^{(0)}_R$ into sub-cubes $\cube_{Rr/\rho}$ of side-length $Rr/\rho$.  We can then split left-hand side of \eqref{transverse-bound} as
$$
( \sum_{\cube_{Rr/\rho}} \| Q( \phi_\kappa, \psi_{\kappa'} ) \|_{L^2(\cube_{Rr/\rho})}^2)^{1/2}$$
(cf. \eqref{p-split}).

Consider the contribution of a single cube $\cube_{Rr/\rho}$.  Intuitively we only expect the portion of $\phi$, $\psi$ which propagates along the backward light cone from $\cube_{Rr/\rho}$ to contribute to this portion.  To make this intuition precise, we will now throw away the tubes $T$, $T'$ which do not intersect $\cube_{Rr/\rho}$.

Consider first the contribution of a single tube $T \in \T(\kappa)$ which does not intersect $\cube_{Rr/\rho}$.  By H\"older, Bessel \eqref{bessel}, and the decay \eqref{local} of $\phi_T$ outside of $T$, we have 
\bas
\| Q( \phi_T, \psi_{\kappa'} ) \|_{L^2(\cube_{Rr/\rho})}
&\lesssim \bigg(\frac{Rr}{\rho}\bigg)^{\frac 12} 
\| \nabla_{x,t} \phi_T \|_{L^\infty_{x,t}}
\| \nabla_{x,t} \psi_{\kappa'} \|_{L^\infty_t L^2_x}\\
&\lesssim \bigg(\frac{Rr}{\rho}\bigg)^{\frac 12} R^{-100n}
E(\psi_{\kappa'})\\
&\lesssim R^{-99n}.
\end{align*}
Summing this over all tubes $T$ and all cubes $\cube_{Rr/\rho}$ we see that this contribution is acceptable.  Thus we may replace $\phi_\kappa$ by $\sum_{T \in \T(\kappa): T \cap \cube_{Rr/\rho} \neq \emptyset} \phi_T$.  Similarly for $\psi_{\kappa'}$ (but we continue to keep the low frequency in $L^\infty_x$ and the high frequency in $L^2_x$).  

We are now ready to apply the hypothesis $B(\alpha)$.  Modulo acceptable errors, we thus have
$$ \| Q( \phi_\kappa, \psi_{\kappa'} ) \|_{L^2(\cube_{Rr/\rho})}
\lesssim (\frac{Rr}{\rho})^\alpha
E(\sum_{T \in \T(\kappa): T \cap \cube_{Rr/\rho} \neq \emptyset} \phi_T)^{1/2} 
E(\sum_{T' \in \T'(\kappa'): T' \cap \cube_{Rr/\rho} \neq \emptyset} \psi_{T'})^{1/2}.$$
Applying \eqref{ortho} we obtain
$$ \| Q( \phi_\kappa, \psi_{\kappa'} ) \|_{L^2(\cube_{Rr/\rho})}
\lesssim (\frac{Rr}{\rho})^\alpha R^{C\eps}
(\sum_{T \in \T(\kappa): T \cap \cube_{Rr/\rho} \neq \emptyset}
\sum_{T' \in \T'(\kappa'): T' \cap \cube_{Rr/\rho} \neq \emptyset}
E(\phi_T) E(\psi_{T'}) )^{1/2}.$$
The expression inside the square root may be thought of as a measure of how much joint energy of $\phi$ and $\psi$ can interact within the cube $\cube_{Rr/\rho}$.  The point is that different cubes $\cube_{Rr/\rho}$ will require disjoint portions of the joint energy because of transversality, and this will allow us to sum.  Indeed, if we insert the above bound into the left-hand side of \eqref{transverse-bound}, we obtain a total bound of
$$ (\frac{Rr}{\rho})^\alpha R^{C\eps}
(\sum_{\cube_{Rr/\rho}} \sum_{T \in \T(\kappa): T \cap \cube_{Rr/\rho} \neq \emptyset}
\sum_{T' \in \T'(\kappa'): T' \cap \cube_{Rr/\rho} \neq \emptyset}
E(\phi_T) E(\psi_{T'}) )^{1/2}.$$
By transversality of $\kappa$ and $\kappa'$ and elementary geometry, we see that for each $T \in \T(\kappa)$, $T' \in \T'(\kappa)$ there are only $O(1)$ cubes $\cube_{Rr/\rho}$ for which the above summation is non-vacuous.  Thus we can bound the previous by
$$ (Rr/\rho)^\alpha R^{C\eps}
(\sum_{T \in \T(\kappa)}
\sum_{T' \in \T'(\kappa')}
E(\phi_T) E(\psi_{T'}) )^{1/2}.$$
Since $\rho \geq \rho_0$, we see that this is acceptable if we set
\be{theta-def}
\rho_0 := R^c r
\end{equation}
for some small absolute constant $c>0$ to be chosen later. 

It remains to control the parallel interactions \eqref{parallel}.  By Bessel's inequality \eqref{bessel} and the same Cauchy-Schwarz trick used to reduce the transverse case to \eqref{transverse-bound}, we may reduce this portion of the estimate to
\be{parallel-bound}
\| Q( \phi_\kappa, \psi_{\kappa'} )
\|_{L^2(\cube^{(0)}_R)} \lesssim R^{(1-c)\alpha} R^{C\eps}
(\sum_{T \in \T(\kappa)} E(\phi_T))^{1/2} (\sum_{T' \in \T'(\kappa')} E(\psi_{T'}))^{1/2} + R^{-10n}
\end{equation}
for all cubes $\kappa \approx \kappa'$ of side-length $l(\kappa) = l(\kappa')= \rho_0$.

Fix $\kappa$, $\kappa'$ and expand
$$
\| Q( \phi_\kappa, \psi_{\kappa'} )
\|_{L^2(\cube^{(0)}_R)} =
\| \sum_{T \in \T(\kappa), T' \in \T'(\kappa') }
Q(\phi_T, \psi_{T'}) \|_{L^2(\cube^{(0)}_R)}.$$
By arguments similar to those in the transverse case, we may use the decay estimate \eqref{local} to insert a characteristic function $\chi_T$ in front of the interaction $Q(\phi_T, \psi_{T'})$ since the error is acceptable.  Similarly we may use \eqref{local-2} to insert a $\chi_{T'}$ in front of $\psi_{T'}$.  Thus we can restrict to those pairs $T$, $T'$ of tubes which intersect.  But in the parallel interaction case we see from the definition \eqref{theta-def} of $\rho_0$ that each $T$ only intersects $O(R^{Cc})$ tubes $T'$.  Also, the tubes $T$ only overlap by a multiplicity of $O(R^{Cc})$.  Thus we may estimate the left-hand side of \eqref{parallel-bound} by
$$
R^{Cc} (\| \sum_{T \in \T(\kappa), T' \in \T'(\kappa'): T \cap T' \neq \emptyset }
Q(\phi_T, \psi_{T'}) \|_{L^2(\cube^{(0)}_R)}^2)^{1/2}.$$
We use the triangle inequality to pull the summations outside the norm. From another Cauchy-Schwartz and application of \eqref{ortho}, \eqref{bessel}, we see that it suffices to show that
$$
R^{Cc} \| Q(\phi_T, \psi_{T'}) \|_{L^2(\cube^{(0)}_R)}
\lesssim 
R^{(1-c)\alpha} E(\phi_T)^{1/2} E(\psi_{T'})^{1/2} + R^{-20n}
$$
for all $T \in \T(\kappa)$ and $T' \in \T'(\kappa')$ such that $T \cap T' \neq \emptyset$.  Indeed we shall show the stronger estimate
\be{parallel-reduce}
\| Q(\phi_T, \psi_{T'}) \|_{L^2(\cube^{(0)}_R)}
\lesssim R^{-\beta} 
E(\phi_T)^{1/2} E(\psi_{T'})^{1/2} + R^{-30n}
\end{equation}
for some absolute constant $\beta$ depending only on $n$, if $c$ was sufficiently small and $r$ was sufficiently close to $\sqrt{R}$.
This shall be deferred to Section \ref{parallel-sec}, after some preliminaries on vector fields.

\section{A digression on vector fields}\label{vector-sec}

In this section we introduce our main tool for obtaining decay and localization estimates, namely the Killing and conformal Killing vector 
fields (see \cite{klainerman1}, \cite{klainerman2} for a more comprehensive discussion).  We begin with an informal discussion.

We are interested in those vector fields $X = X^\alpha \partial_\alpha$ which preserve the space of free waves.  In addition to the constant co-efficient vector fields $(\partial_1, \ldots, \partial_n, \partial_t)$, one also has the Poincare group of vector fields
$$
\left.
\begin{array}{rllr}
S&:=&t\pr_t+ r \pr_r & \hbox{(Scaling)}\\
L_i&:=&x_i\pr_t+t\pr_i \hbox{ for } 1 \leq i \leq n& \hbox{(Lorentz boosts)}\\
\Om_{i,j}&:=&x_i\pr_j-x_j\pr_i. \hbox{ for } 1 \leq i < j \leq n &
\hbox{(Angular rotations)} 
\end{array}
\right.
$$
We use 
$$\Gamma = \Gamma^{(0,0)} := (S, L_1, \ldots, L_n, \Omega_{1,2}, \ldots, \Omega_{n-1,n})$$
to denote the entire collection of Poincare vector fields, thus
$\Gamma$ is a vector-valued first order differentiation operator.  We can then define $\Gamma^k$ for any $k \geq 1$ as a $k$-vector-valued $k^{th}$ order differentiation operator.

The Poincare vector fields $\Gamma = \Gamma^{(0,0)}$ are centered at the origin, but we may also translate them in spacetime to be centered at an arbitrary point $(x_0, t_0)$, creating a new set of conformal Killing vector fields $\Gamma^{(x_0,t_0)}$.  This will turn out to be important for tube localization.

If $\phi$ is a free wave and $D$ is a collection of conformal Killing vector fields (e.g. $D = \Gamma$), then $D^k \phi$ is also a free wave, and in particular we have the energy identity
$$ \| \nabla_{x,t} D^k \phi(t) \|_2 = \| \nabla_{x,t} D^k \phi(0) \|_2.$$
To exploit this we seek to use vector fields $D$ which are ``large'' at time t and ``small'' at time zero.  In practice, this means that we try to ensure $\phi(0)$ is localized in physical space (in order to make the weights associated with the $X_j$ small at time zero), and at time $t$ we interest ourselves with the regions where the vector fields $D$ are large and span the tangent space.  In such regions we will obtain extremely good decay estimates\footnote{By taking $k$ arbitrarily large, we will in fact obtain arbitrarily fast decay; we will be able to justify taking these many derivatives because we will always be working in the frequency-localized setting.  To some extent this is still true in such apparently rough situations as quasilinear wave equations, because of the frequency localization.  In this paper we shall differentiate many times (e.g. $100n$ times) but one can get away with far less differentiation by being more careful with the error terms, etc.}.  We observe the convenient fact that if $\phi$ has frequency $\mu$ and $D$ is any collection of constant co-efficient or Poincare vector fields, then $D^k \phi$ also has frequency $\mu$.

To localize a wave to a cone centered at the origin $(0,0)$ it is enough to use the Poincare vector fields $\Gamma^{(0,0)}$ (in fact one can use only the scaling vector field $S$
to do this).  Similarly for cones centered around other points $(x_0,t_0)$ in spacetime.  To localize instead to a tube, say one with velocity $e_1$, one shall need some additional (constant coefficient) vector fields, notably $\partial_2, \ldots, \partial_n$ and $\partial_t + \partial_1$.

We now turn to the details.  
Let $\psi$ be a free wave with some frequency $\mu$.  Then $\Gamma^k \psi$ is also a free wave with frequency $\mu$, and the energy identity yields
$$
\| \nabla_{x,t} \Gamma^k \psi(t) \|_2 \lesssim \| \nabla_{x,t} \Gamma^k \psi(0) \|_2$$
for all $k \geq 0$ (here and in the sequel all implicit constants may depend on $k$).  On the other hand, if we expand $\Gamma^k$ out using the Leibnitz product rule, and use the equation $\Box \psi = 0$ to convert any repeated time derivative to a repeated space derivative, we obtain the pointwise estimate
$$
 |\nabla_{x,t} \Gamma^k \psi(x,0)| \lesssim \sum_{j=0}^k |x|^j |\nabla_x^j \nabla_{x,t} \psi(x,0)|.
$$
Combining these bounds with \eqref{pointwise} we obtain the weighted energy estimate
\be{weighted}
\| \nabla_{x,t} \Gamma^k \psi(t) \|_2 \lesssim \| (1 + \mu |x|)^{k}
\nabla_{x,t} \psi(x,0) \|_{L^2_x}.
\end{equation}

As a first attempt to exploit this differentiability in $\Gamma$ we observe the identities
\bas
\partial_t &= \frac{t}{t^2 - r^2} S - \sum_{j=1}^n \frac{x_j}{t^2 - r^2} L_j\\
\partial_i &= \frac{1}{t} L_i - \frac{x_i}{t} \partial_t \\
&= \frac{1}{t} L_i - \frac{x_i}{t^2 - r^2} S + \sum_{j=1}^n \frac{x_i x_j}{t(t^2 - r^2)} L_j
\end{align*}
where we use polar co-ordinates $x = r\omega$
As a caricature, we thus have something like
$$ \nabla_{x,t} \lq\lq\sim{\rq}{\rq} \frac{x^2 + t^2}{t(t^2 - r^2)} \Gamma \lq\lq\lesssim{\rq}{\rq} \frac{1}{||t|-|x||} \Gamma.$$
If we combine this with repeated applications of the Leibnitz product rule, we conclude the pointwise estimate
\be{away-from-cone}
|\nabla_{x,t}^k \psi(x,t)| \lesssim \frac{1}{||t|-|x||^k} \sum_{j=0}^k   | \Gamma^j \psi(x,t) |.
\end{equation}
This is estimate is good away from the light cone, but is not completely satisfactory in the region $|x|/2 < |t| < 2|x|$.  In this case we introduce the null frame
$$
\left.\begin{array}{rlr}
E_+ &:= \sgn(t) \pr_t + \pr_r & \hbox{good null direction}\\
E_- &:= \sgn(t) \pr_t - \pr_r & \hbox{bad null direction}\\
A_i &:= \pr_i-\omega_i \pr_r & \hbox{angular derivatives}.
\end{array}\right.
$$
We write $\nabb := (A_1, \ldots, A_n)$ for the angular derivatives, and
$\dabb := (E_+, \nabb)$ for the good derivatives.

These vector fields are not quite conformal Killing, but we can write them in terms of the Poincare vector fields:
\bas
(t+r) E_+ &= S + \sum_{i=1}^n \omega_i L_i \\
(t-r) E_- &= S - \sum_{i=1}^n \omega_i L_i \\
t A_i &= L_i - \frac{x_i}{r} \sum_{j=1}^n \omega_i L_j.
\end{align*}
As a caricature, we thus have
$$ \dabb ``\lesssim{\rq}{\rq} \frac{1}{|t|} \Gamma; \quad E_- ``\lesssim{\rq}{\rq} \frac{1}{||t|-|x||} \Gamma.$$
By repeated applications of the Leibnitz rule, we can thus improve \eqref{away-from-cone} to 
\be{k-decay}
|\dabb^j \nabla_{x,t}^k \psi(t,x)| \lesssim \frac{1}{|t|^j \big||t|-|x|\big|^k} \sum_{m=0}^{j+k} |\Ga^m \psi(t,x)|
\end{equation}
near the light cone for all $j, k \geq 0$.  Informally, this means that if the $\Gamma$ vector fields are under control, then good derivatives gain a $|t|$ while bad derivatives only gain a $||t|-|x||$.

We now specialize the above estimates.  Fix $R \gtrsim 1$, $\mu \gtrsim 1$, and suppose that $\psi$ is a wave of frequency $\mu$ such that $E(\psi) \lesssim 1$.  We also assume that $\psi$ is localized around the origin at time zero.  More precisely, we assume there is a radius $1 \lesssim \Lambda \ll R$ such that we have the estimate
\be{good-decay}
(\int_{|x| \gg \Lambda} |x|^{200n} |\nabla_{x,t} \psi(x,0)|^2\ dx)^{1/2} \lesssim 1. 
\end{equation}
From \eqref{weighted} we thus obtain 
$$ \| \nabla_{x,t} \Gamma^k \psi(0) \|_2 \lesssim \Lambda^{k} \mu^{k} E(\psi)^{1/2} + \Lambda^{-100n + k} \mu^{k}$$
for all $0 \leq k \leq 100n$.  Since $\Gamma^k \psi$ has frequency $\mu$, we can use the reproducing formula \eqref{reproducing} to then obtain 
\be{decay-vector}
\| \Gamma^k \psi(t) \|_2 \lesssim \Lambda^k \mu^{k-1} E(\psi)^{1/2} + \Lambda^{-100n + k} \mu^{k-1}.
\end{equation}

We now specialize to the case when $\psi$ has frequency $\mu = 1$, and time is restricted to the region $R \leq t \leq 2R$.  To emphasize the frequency 1 we shall write $\phi$ instead of $\psi$.

First suppose that we are in the region $r \leq |t|/2$ or $r \geq 2|t|$; in other words, we are away from the cone.  By \eqref{away-from-cone}, \eqref{decay-vector} we see that
\be{rapid}
|\nabla_{x,t}^k \phi(t,x)| \lesssim \Lambda^k R^{-k}  E(\phi)^{1/2} + \Lambda^{-100n}.
\end{equation}
Now suppose that we are in the region $|t|/2 < r < 2|t|$, so that $r \sim R$. From \eqref{k-decay} we have
\be{star}
|\nabb^j \nabla_{x,t}^k \phi(t,x)| \lesssim R^{-j} \sum_{m=0}^j |\Gamma^{m} \nabla_{x,t}^k \phi(t,x)|.
\end{equation}
We now use the angular Sobolev inequality 
$$ \| u \|_{L^\infty(S^{n-1})} \lesssim \| u \|_{L^1(S^{n-1})}
+ \| \nabb^{n-1} u \|_{L^1(S^{n-1})},$$
which is easily verified by $n-1$ applications of the fundamental theorem of Calculus.  Applying this inequality to $u := |\nabla_{x,t}^k \phi(t,x)|^2$, we obtain
$$
|\nabla_{x,t}^k \phi(t,x)|^2 \lesssim \int_{S^{n-1}}
\sum_{i+j \leq n-1}
|\nabb^i \nabla_{x,t}^k \phi(t,r\omega')|
|\nabb^j \nabla_{x,t}^k \phi(t,r\omega')|\ d\omega'.$$
By \eqref{star} we thus have
$$
|\nabla_{x,t}^k \phi(t,x)|^2 \lesssim (\frac{\Lambda}{R})^{n-1} \int_{S^{n-1}}
r^{n-1} (\sum_{0 \leq i \leq n-1}
\Lambda^{-i} |\Gamma^i \nabla_{x,t}^k \psi(t,r\omega')|)^2\ d\omega'.$$
This is a square of the $L^2$ norm on a sphere of radius $r$.
We then use \eqref{pointwise} to average out this sphere into what is essentially an annulus of thickness $O(1)$ around the sphere of radius $r$.  This implies that
$$
|\nabla_{x,t}^k \phi(t,x)|^2 \lesssim (\frac{\Lambda}{R})^{n-1}
\sum_{0 \leq i \leq n-1} \Lambda^{-i} E(\Gamma^i \nabla_{x,t}^{k-1} \phi)^{1/2}.$$
Applying \eqref{decay-vector} we thus obtain the pointwise estimate
$$ 
|\nabla_{x,t}^k \phi(t,x)| \lesssim  (\frac{\Lambda}{R})^{\frac{n-1}{2}} E(\phi)^{1/2}
+ \Lambda^{-90n}.$$
If we combine this with \eqref{rapid} we obtain the decay estimate
\be{decay-final}
\| \nabla_{x,t}^k \phi(t) \|_\infty \lesssim (\frac{\Lambda}{R})^{\frac{n-1}{2}} E(\phi)^{1/2} + \Lambda^{-90n}
\end{equation}
for $k=100n$ (say).  But since $\psi$ has frequency $1$, this estimate in fact holds for all $k$ by the reproducing formula \eqref{reproducing}.  Finally, we remark that we may replace $\phi$ by $\Gamma \phi$ on both sides without affecting the estimate.  In particular we have
\be{decay-gamma}
\| \Gamma \phi(t) \|_\infty \lesssim (\frac{\Lambda}{R})^{\frac{n-1}{2}} E(\Gamma \phi)^{1/2} + \Lambda^{-90n} 
\lesssim (\frac{\Lambda}{R})^{\frac{n-1}{2}} \Lambda E(\phi)^{1/2} + \Lambda^{-90n}.
\end{equation}

\section{Parallel interactions}\label{parallel-sec}

Fix $T$, $T'$.  To complete the proof of Theorem \ref{null-1} begin in Section \ref{induction-sec}, we need to prove \eqref{parallel-reduce}.  (We also need to prove Lemma \ref{tube-decomp}, but this will be done in later sections).

In earlier work such as \cite{tv:cone2}, \cite{wolff:cone}, \cite{tao:cone} the proof of \eqref{parallel-reduce} was accomplished by using Lorentz transforms to reduce to the $R=1$ case.  However we do not wish to do this here as this method seems to rely heavily on the special properties of Minkowski space.  Instead we shall use commutation with vector fields.

In proving \eqref{parallel-reduce} we will actually prove a somewhat stronger result:

\begin{proposition}\label{par-interact} Let $\sqrt{R} \ll \Lambda \ll R$, and $\phi$, $\psi$ be free waves of frequency $1$ and $\mu$ respectively, where $\mu \gtrsim 1$.  Suppose we have the energy normalization
$$ E(\phi), E(\psi) \lesssim 1$$
and the localization properties
$$
(\int_{|x-x_0| \gg \Lambda} |x-x_0|^{200n} |\nabla_{x,t} \phi(x,t_0)|^2\ dx)^{1/2} \lesssim 1$$
and
$$
(\int_{|x-x_0| \gg \Lambda} |x-x_0|^{200n} |\nabla_{x,t} \psi(x,t_0)|^2\ dx)^{1/2} \lesssim 1
$$ 
for some $x_0 \in \R^n$ and $0 \leq t_0 \leq R/2$.  Then we have the fixed-time estimate
\be{fixed}
\| Q(\phi,\psi)(t) \|_{L^2(\R^n)} \lesssim (\frac{\Lambda}{R})^{(n+1)/2} E(\phi)^{1/2} E(\psi)^{1/2} + R^{-40n}
\end{equation}
for all $R \leq t \leq 2R$.
\end{proposition}

\begin{proof}
By translation invariance we may set $x_0 = 0$.  Also by a time translation and some redefinition of $R$ if necessary we may set $t_0 = 0$.

Fix $\phi$, $\psi$, $R$, $\mu$.  In the region away from the light cone, i.e. $|x| \leq t/2$ or $|x| \geq 2t$, we can control $\phi$ using \eqref{rapid} and $\psi$ using energy estimates to easily obtain \eqref{fixed}.  Now suppose we are near the light cone, so that $|x| \sim t \sim R$. 

We shall write the null form $Q$ in terms of the null frame $A_i$, $E_+$, $E_-$ developed in the previous section.  More precisely, we claim the pointwise estimate
$$ |Q(\phi,\psi)| \lesssim |\dabb\phi| |\nabla_{x,t} \psi| + |\nabla_{x,t} \phi| |\dabb\psi|.$$
To see this, we substitute in the identities
$\partial_i = A_i + \frac{x_i}{2r} E_+ - \frac{x_i}{2r} E_-$
and $\partial_t = \frac{1}{2}E_+ + \frac{1}{2} E_-$, and observe that in each of the three null forms $Q_0$, $Q_i$, $Q_{ij}$ the $E_- \phi E_- \psi$ term (which is the only potentially bad term) always vanishes\footnote{In principle the $Q_0$ null form should be somewhat better than $Q_i$ and $Q_{ij}$, as can be seen by inspection of the bilinear symbol (see e.g. \cite{klainerman:foschi}), but this cannot be detected by this simple null frame analysis.}.

From the previous, \eqref{k-decay} and H\"older we have
$$ \| Q(\phi,\psi)(t) \|_2 \lesssim \frac{1}{R} 
(\| \Gamma \phi(t)\|_\infty \| \nabla_{x,t} \psi(t)\|_2 +
 \|\nabla_{x,t} \phi\|_\infty \|\Gamma \psi(t)\|_2).$$
The claim now follows from estimating $\phi$ using decay estimates \eqref{decay-final}, \eqref{decay-gamma}, and $\psi$ using the standard energy estimate and \eqref{decay-vector}.
\end{proof}

We now apply this Proposition with $\phi$, $\psi$ replaced by $\phi_T$ and $\psi_{T'}$ respectively, with $t_0$ set to $R/2$, and $x_0$ set to the center of $\kappa$.  Recall that $\rho_0= rR^c$ is the size of $\kappa$. By the Proposition and an $L^2$ norm in time we have
$$ \| Q(\phi_T,\psi_{T'}) \|_{L^2(\R^n \times [R, 2R]))} \lesssim R^{1/2} R^{Cc} (\frac{r}{R})^{(n+1)/2} E(\phi_T)^{1/2} E(\psi_{T'})^{1/2} + R^{-40n}.$$
If $c$ is sufficiently small, and $r$ is sufficiently close to $\sqrt{R}$, the claim follows.  (Indeed from the argument we see that we just need $r$ to be substantially less than $R^{n/(n+1)}$).

\section{First proof of Lemma \ref{tube-decomp}: Phase space preparation of the initial data.}\label{first-proof}

Fix $\phi$, $\mu$, $R$, $r$.  It suffices to prove the claim under the assumption on the initial data that $\phi_t(0) = \pm i\sqrt{-\Delta} \phi(0)$ for some sign $\pm$; this is because general initial data can be split into two pieces of the above form via the formula
$$ \phi[0] = \sum_\pm
(\frac{\phi(0) \pm i\sqrt{-\Delta}^{-1} \phi_t(0)}{2},
\pm \sqrt{-\Delta} \frac{\phi(0) \pm i\sqrt{-\Delta}^{-1} \phi_t(0)}{2})$$
Note that this decomposition is orthogonal with respect to the energy inner product $\langle \cdot,\cdot \rangle_e$, and does not affect the property of $\phi$ having frequency $\mu$.

We shall assume that $\phi_t(0) = + i\sqrt{-\Delta} \phi(0)$; the argument for the other sign is similar and is left to the reader.  It shall be convenient to define the length scales $s := R^{-\eps/2} r$ and $\Lambda := R^{-\eps/4} r$.  

Since the $\phi_T$ are free waves, to construct $\phi_T$ it will suffice to construct their initial data $\phi_T[0]$.  By standard arguments we may partition $1 = \sum_{x_0 \in \Sigma} \chi_{x_0}$ on $\R^n$, where $\Sigma \subset \R^n$ is a lattice of spacing $\sim s$, and $\chi_{x_0}$ is a bounded function with Fourier transform supported in the ball $|\xi| \ll \mu$ and which satisfies the estimates
\be{nabla-bound}
|\chi_{x_0}(x)| \lesssim R^{-M} \mu^{-M} \hbox{ whenever } |x-x_0| \gg s
\end{equation}
for any $M > 0$; the implicit constants may depend on $M$.
For instance, one can set $\chi_{x_0} := P_{\mu/C} \chi_{\cube(x_0)}$, where the $\cube(x_0)$ form a partition of $\R^n$ into cubes of side-length $\sim s$.

We also partition $1 = \sum_{\omega \in \Omega} \eta_\omega$ on $S^{n-1}$, where $\Omega$ is a $s/R$-separated set of directions and $\eta_\omega$ is a bump function adapted to the disk of radius $O(s/R)$ centered at $\omega$.

We define the sector multipliers $P_{\mu,\omega}$ via the Fourier transform as
$$ \widehat{P_{\mu,\omega} f}(\xi) := \eta(\frac{\xi}{|\xi|}) a(|\xi|) \hat f(\xi)$$
where $a$ is a bump function which equals one on a suitable annulus $|\xi| \sim \mu$ and is adapted to a dilate of this annulus.  Thus $P_{\mu,\omega}$ is a Fourier projection to a tubular region of length $\mu$ and width $s\mu/R$, oriented in the direction $\omega$.  Observe from construction that $s \cdot \mu$ and $s\mu/R \cdot s$ are both $\gtrsim R^\eps \mu$.  By the uncertainty principle we thus see that the convolution kernel $K_\omega$ of $P_{\mu,\omega}$ is integrable and obeys the estimates
\be{nabla-K-bound}
|\nabla K_\omega(x)| \lesssim R^{-M} \mu^{-M} \hbox{ whenever } |x| \gtrsim s
\end{equation}
for any $M$.

Our collection of tubes $\T$ will be defined by
$$ \T := \{ T(R,r,x_0,\omega): x_0 \in \Sigma; \omega \in \Omega \};$$
note that the claimed cardinality bound on the intersection of $\T$ with any $R$-cube is clear.  

Next, our $\phi_T$ will be defined via their initial data as
$$ \phi_{T(R,r,x_0,\omega)}(0) := P_{\mu,\omega}( \chi_{x_0} \phi(0) )$$
and
$$ \partial_t \phi_T(0) := i\sqrt{-\Delta} \phi_T(0).$$
We will not need an error term for this argument and will set $\phi_{error} := 0$.

By construction it is clear that the $\phi_T$ are free waves with frequency $\sim \mu$, and that \eqref{decomposition} holds.  

Observe that the physical space multipliers $\chi_{x_0}$ are almost orthogonal, as are the Fourier space multipliers $P_{\mu,\omega}$.  Thus \eqref{bessel} holds.  Now we show \eqref{ortho}. Since the time derivative of $\phi_T$ is controlled by the spatial derivatives, it suffices to show that
$$ \| \sum_{(\omega, x_0) \in X} \nabla_x P_{\mu,\omega}( \chi_{x_0} \phi(0) ) \|_2
\lesssim R^{C\eps}
(\sum_{(\omega, x_0) \in X} \| \nabla_x P_{\mu,\omega}( \chi_{x_0} \phi(0) ) \|_2^2)^{1/2} + R^{-100n}$$
for all subsets $X$ of $\Omega \times \Sigma $.

Observe from Plancherel's theorem that the $\nabla_x P_{\mu,\omega}$ multipliers are essentially orthogonal as $\omega$ varies.  Thus it suffices to show that
$$ \| \sum_{x_0 \in Y} \nabla_x P_{\mu,\omega}( \chi_{x_0} \phi(0) ) \|_2
\lesssim R^{C\eps}
(\sum_{x_0 \in Y} \| \nabla_x P_{\mu,\omega}( \chi_{x_0} \phi(0) ) \|_2^2)^{1/2} + R^{-150n}$$
for all $\omega \in \Omega$ and $Y \subset \Sigma$.

Fix $\omega$, $Y$.  From \eqref{nabla-bound}, \eqref{nabla-K-bound} and the assumption $E(\phi) \lesssim 1$ we see that
\be{joy}
|\nabla_x P_{\mu,\omega}( \chi_{x_0} \phi(0) ) (x)| \lesssim (R + |x-x_0|)^{-M} \mu^{-M} \hbox{ whenever } |x - x_0| \gtrsim \Lambda
\end{equation}
for any $M$.
Thus one can restrict $\nabla_x P_{\mu,\omega}( \chi_{x_0} \phi(0) ) (x)$ to a ball of radius $O(\Lambda)$ around $x_0$ by paying a total price of $O(R^{-150 n})$.  Since these balls only overlap by about $O(R^{C\eps})$ or so, the claim follows.  This proves \eqref{ortho}.

It remains to show that each $\phi_T$ is localized to $T$.  Fix $T$; without loss of generality we may take $x_0 = 0$ and $\omega = e_1$.  We have to show \eqref{local}, \eqref{local-2}.  We only show the $L^2$ bound \eqref{local-2}; the $L^\infty$ bound \eqref{local} essentially follows from Sobolev (replacing $r$ by $r/2$ if necessary).  We thus need to show
\be{local-decay}
(\int_{|x-te_1| \geq r} | \nabla_{x,t} \phi_T(x,t) |^2\ dx)^{1/2} \lesssim (R + |x - t e_1|)^{-100n}
\end{equation}
whenever $|t| \sim R$.  

We shall achieve this using vector fields by means of the following lemma, which will also be used in our second proof in the next section.

\begin{lemma}\label{local-lemma}  Let $R^{1/2+\eps} \lesssim r \ll R$, and let $\Lambda = R^{-\eps/4} r$.  Let $\phi_T$ be a free wave of some frequency $\mu \gtrsim 1$ such that
\be{point-bound}
| \nabla_{x,t}^k \phi_T(x,0) | \lesssim (R + |x|)^{-M} \mu^{-M} \hbox{ whenever } |x| \gtrsim \Lambda
\end{equation}
for any $k \geq 1$ and any $M$, with the implicit constant depending on $k$ and $M$ of course.  Also, if $\Sigma$ denotes the vector fields
$$ \Sigma := \{ \partial_t + \partial_1, \partial_2, \ldots, \partial_n \}$$
then we assume the estimates
\be{sigma-bound}
E(\Sigma^k \phi_T) \lesssim R^{100n} (\frac{\Lambda \mu}{R})^k
\end{equation}
for all $k \geq 0$.  Then \eqref{local-decay} holds for all $|t| \sim R$.
\end{lemma}

\begin{proof}
Since $\phi_T$ has frequency $\mu$, it will suffice to show that
\be{chip}
(\int_{|x-te_1| \geq r} | \nabla_{x,t}^{k+1} \phi_T(x,t) |^2\ dx)^{1/2}
 \lesssim (R + |x - t e_1|)^{-100n} \mu^k
\end{equation}
for some suitable $k$ which we can choose at will.

Fix $t$.  From \eqref{point-bound} and \eqref{sigma-bound} we have the regularity estimate 
\be{shoop}
E(\Gamma^k \Sigma^j \phi_T)^{1/2} \lesssim R^{100n} (\Lambda \mu)^k (\frac{\Lambda \mu}{R})^j + R^{-M} \mu^{-M}
\end{equation}
for any $k, j \geq 0$.

We now divide the region $|x-te_1| \geq r$ into several sub-regions and prove the localization estimate \eqref{chip} separately for each sub-region.

{\bf Case 1. $||t|-|x|| \gtrsim r$}

In this case we are away from the light cone and one can use the fact that the vector fields $\Gamma$ span.

By \eqref{k-decay} we have
$$ |\nabla_{x,t}^{k+1} \phi_T(x,t)| \lesssim r^{-k} |\Gamma^k \nabla_{x,t} \phi_T(x,t)|$$
and hence by the regularity estimate \eqref{shoop}
$$ (\int_{x: ||t|-|x|| \gtrsim r} |\nabla_{x,t}^{k+1} \phi_T(x,t)|^2\ dx)^{1/2}
\lesssim R^{100n} (\Lambda \mu/r)^k + R^{-M} \mu^{-M}.$$
Since $r = R^{\eps/4} \Lambda$, we thus obtain \eqref{chip} for this region by setting $k$ sufficiently large.

{\bf Case 2. $||t|-|x|| \ll r$, and $x_1 > |t|/2$.}

Without loss of generality we may take $|x_2| \geq |x_3|, \ldots, |x_n|$, thus $|x_2| \gtrsim r$.  In particular $|\omega_2| \gtrsim r/R$.  The idea is now to use the fact that $\Gamma$ together with $\partial_2$ spans.

We now use the identity
$$ E_- = E_+ + \frac{2}{\omega_2} A_2 - \frac{2}{\omega_2} \partial_2;$$
since $E_-$ and $\dabb$ form a null frame,
we thus informally have
$$ \nabla_{x,t} \lq\lq\lesssim{\rq}{\rq} (1 + \frac{1}{\omega_2}) (\dabb, \Sigma).$$
Since $|\omega_2| \gtrsim r/R$, we thus obtain after repeated applications of the Leibnitz rule
$$ |\nabla^{k+1}_{x,t} \phi_T(x,t)| \lesssim \sum_{i+j \leq k} R^{i+j} r^{-k} |\nabla_{x,t} \dabb^i \Sigma^j \phi_T(x,t)| $$
for any $k$.  Applying the regularity estimate \eqref{shoop} and substituting $\Lambda = R^{-\eps/4} r$, we thus obtain
$$ (\int_{|x_2| \gtrsim r} |\nabla^{k+1}_{x,t} \phi_T(x,t)|^2\ dx)^{1/2} \lesssim R^{100n} \mu^k R^{-\eps k/2}
E(\phi_T)^{1/2} + R^{-M} \mu^{-M}.$$
The claim follows by setting $k$ sufficiently large.

{\bf Case 3. $||t|-|x|| \ll r$, and $x_1 \leq |t|/2$.}

In this case we use the fact that $\dabb$ together with $\partial_t + \partial_1$ span.  Specifically, we use the identity
$$ E_- = \frac{-1-\omega_1}{1-\omega_1} E_+
+ \frac{2}{1-\omega_1} (\partial_t + \partial_1) - \frac{2}{1-\omega_1} A_1.$$
Since $|1-\omega_1| \sim 1$ by hypothesis, we thus informally have
$$ \nabla_{x,t} ``\lesssim'' (\dabb, \Sigma).$$
One then argues as in Case 2.
\end{proof}

Fix $t$.  We now apply the Lemma just proven to obtain the localization estimate \eqref{local-decay}.
By repeating the derivation of the initial data localization estimate \eqref{joy} we can obtain the pointwise estimates \eqref{point-bound}.  We now verify the regularity estimate \eqref{sigma-bound}.

For the derivatives $\partial_2, \ldots, \partial_n$ this is clear from the Fourier support of $\phi_T$.  Now consider the powers of $\partial_t + \partial_1$.
By repeated application of the equation $\Box \phi_T = 0$ we see that
$$ \partial_t^k \phi_T = (i\sqrt{-\Delta})^k \phi_T,$$
and hence
$$ (\partial_t + \partial_1)^k \phi_T = (i\sqrt{-\Delta} - \partial_1)^k \phi_T.$$
On the other hand, the Fourier transform of $\nabla_{x,t}^k \phi_T$ is supported in the region $\xi_1 \sim \mu$, $\xi = \xi_1 e_1 + O(\mu s/R)$.  In particular we have the crude estimate $i|\xi| - i\xi_1 = O(\mu s/R)$, and the claim follows.  This concludes the proof of the localization estimate \eqref{local-decay}, and Lemma \ref{tube-decomp} follows. 

\section{Second proof of Lemma \ref{tube-decomp}: Spatial truncation at two distinct times.}\label{second-proof}

In the previous argument, tube localization was achieved by localizing $\phi[0]$ in physical space and then applying a Fourier projection.  Although this method works, it is somewhat unsatisfactory as it (implicitly) uses the fact that the Fourier support of a free wave is preserved by the time evolution.  This type of fact does not seem very stable in rough metrics, and so we would like to use a tube decomposition method which does not require as much usage of Fourier space.

The idea is to characterize a tube by its intersection with two time slices $t=0$ and $t=3R$ (cf. the strategy used to attack the Kakeya problem in \cite{borg:high-dim}).  We shall write $\tau$ for $3R$.

In this section it will be convenient to take $s := R^{1/2 + \eps/2}$ and $\Lambda := R^{1/2 + 3\eps/4}$.
With this value of $s$, we may partition $1 = \sum_{x_0 \in \Sigma} \chi_{x_0}$ on $\R^n$ as in the previous section.

Let $U(t)$ denote the free evolution operator on free waves, thus $U(t) \phi[0] := \phi[t]$.  Then we may use linearity to split
$$ \phi[t] = \sum_{x_0, x_\tau \in \Sigma}  U(t)( \chi_{x_0} U(-\tau) (\chi_{x_\tau} \phi[\tau])).$$
In other words, for each $x_0$, $x_\tau$, we localize $\phi$ to an $s$-neighbourhood of $x_\tau$ at time $\tau$, evolve back by the free wave equation, localize to an $s$-neighbourhood of $x_0$ at time $0$, and use this as initial data.  As we shall see (using the fact that $s > \sqrt{R}$), the second localization will not significantly disturb the localization of the first.  

We divide the summation into two terms, depending on whether we have a null ray (in the sense that $|x_\tau - x_0| = \tau + O(\Lambda)$) or not.

First consider the contribution of the case $||x_\tau| - |x_0|| \neq \tau + O(\Lambda)$, so that the two points do not form a null ray.  This case is small and will be absorbed into the error term $\phi_{error}$.  To prepare for this, we fix $x_\tau$, and define $\psi = \psi_{x_\tau}$ by $\psi[t] := U(t-\tau) (\chi_{x_\tau} \phi[\tau])$.  Thus $\psi$ is a free wave of frequency $\mu$ which is essentially localized to the disk $\{(x,\tau): x = x_\tau + O(s)\}$ at time $\tau$.  To exploit this we recall that $\Gamma^{(x_\tau,\tau)}$ denotes the Poincare vector fields centered at the point $(x_\tau,\tau)$ in spacetime.  From the reproducing formula \eqref{reproducing} and several applications of the Leibnitz rule we have 
$$ E((\Gamma^{(x_\tau,\tau)})^k \psi) \lesssim (s\mu)^{2k} \int (1 + \frac{|x - x_\tau|}{\Lambda})^{-100n} |\nabla_{x,t} \phi(x,\tau)|^2\ dx$$
for any $k$, where of course the implicit constant depends on $k$.  By the pointwise estimates \eqref{k-decay} (adapted to $\Gamma^{(x_\tau,\tau)}$ rather than $\Gamma^{(0,0)}$) we can thus bound 
$$ (\int_{||x|-\tau| \gtrsim \Lambda} |\nabla_{x,t}^{k+1} \psi(x,0)|^2\ dx)^{1/2}
\lesssim \Lambda^{-k} (s\mu)^{k} 
(\int (1 + \frac{|x - x_\tau|}{\Lambda})^{-100n} |\nabla_{x,t} \phi(x,\tau)|^2\ dx)^{1/2}.$$
Since $\psi$ has frequency $\mu$, we thus see from the pointwise bounds \eqref{pointwise} that
$$ (\int_{||x|-\tau| \gtrsim \Lambda} |\nabla_{x,t}\psi(x,0)|^2\ dx)^{1/2}
\lesssim (s/\Lambda)^{k} 
(\int (1 + \frac{|x - x_\tau|}{\Lambda})^{-100n} |\nabla_{x,t} \phi(x,\tau)|^2\ dx)^{1/2}.$$
By making $k$ sufficiently large we may set $(s/\Lambda)^k \lesssim R^{-100n}$.  From this estimate it is easy to show that
$$ E( \sum_{x_0,x_\tau \in \Sigma: ||x_\tau| - |x_0|| \neq \tau + O(\Lambda)}
\chi_{x_0} U(-\tau) (\chi_{x_\tau} \phi[\tau]) ) \lesssim R^{-50n}$$
which is acceptable.

It thus remains to consider the terms when $||x_\tau| - |x_0|| = \tau + O(\Lambda)$.  In this case we can find a tube $T = T(x_0,x_\tau)$ which is a distance $O(\Lambda)$ from $(x_0,0)$ and $(x_\tau,\tau)$, and we define $\phi_T(t) = U(t)( \chi_{x_0} U(-\tau) (\chi_{x_\tau} \phi[\tau]))$ accordingly.  If we let $\T$ denote the collection of all such tubes (which we can of course make disjoint) then \eqref{decomposition} holds.  To show \eqref{bessel}, we observe from spatial orthogonality and the energy identity that
\bas 
\sum_{T \in \T} E(\phi_T)
&\lesssim
\sum_{x_0, x_T \in \Sigma}
E(U(t)( \chi_{x_0} U(-\tau) (\chi_{x_\tau} \phi[\tau]))) \\
&=
\sum_{x_0, x_\tau \in \Sigma}
E(\chi_{x_0} U(-\tau) (\chi_{x_\tau} \phi[\tau])) \\
&\lesssim
\sum_{x_\tau \in \Sigma}
E(U(-\tau) (\chi_{x_\tau} \phi[\tau])) \\
&=
\sum_{x_\tau \in \Sigma}
E(\chi_{x_\tau} \phi[\tau]) \\
&\lesssim 
E(\phi[\tau]) = E(\phi)
\end{align*}
as desired.  

Now we show the $L^2$ localization estimate \eqref{local-2}; the $L^\infty$ localization estimate \eqref{local} is similar and follows from \eqref{local-2} (with $r$ replaced by $r/2$) and the pointwise estimate \eqref{pointwise}.  

Without loss of generality we shall assume that $x_0 = 0$ and $x_\tau = \tau e_1 + O(\Lambda)$.  In this case we need to prove the localization estimate \eqref{local-decay}.  By Lemma \ref{local-lemma} it suffices to verify the initial data localization estimate \eqref{point-bound} and the regularity estimate \eqref{sigma-bound}.

The claim \eqref{point-bound} follows easily from the decay of cutoff $\chi_{x_0}$ and crude energy estimates on $\phi_T$ (the Bessel inequality\eqref{bessel} will do).  Now we verify \eqref{sigma-bound}.  Once again we set $\psi = \psi_{x_\tau}$ by $\psi[t] := U(t-\tau) (\chi_{x_\tau} \phi[\tau])$.   If derivatives $\Sigma$ land on the cutoff $\chi_{x_0}$ then one gains a power of $s$ for each derivative (note that we can use the equation $\Box \phi_T=0$ to convert any repeated time derivative to a repeated space derivative).  Thus it suffices to show that
\be{reg}
E(\chi_{x_0} \Sigma^k \psi[0]) \lesssim R^{100n} (\frac{\Lambda \mu}{R})^k.
\end{equation}
We may replace $\chi_{x_0}$ with the characteristic function of the disk $D := \{ (x,0): |x-x_0| \lesssim \Lambda \}$ since the rapid decay of $\chi_{x_0}$ ensures that the contribution outside this disk is acceptable.

This disk is close to the light cone centered at $(x_\tau,\tau)$.  We will thus use the good derivatives $\dabb^{(x_\tau,\tau)}$, being the vector fields $\dabb$ translated to be centered at $(x_\tau,\tau)$ and adapted to the backward time flow.

We observe the identities
\bas
 \partial_t + \partial_1 &= A_1^{(x_\tau,\tau)} - \frac{1}{2}(1 - \omega_1^{(x_\tau,\tau)}) E_+^{(x_\tau,\tau)} 
- \frac{1}{2}(1 + \omega_1^{(x_\tau,\tau)}) E_-^{(x_\tau,\tau)}\\
\partial_i &= A_i +\frac{1}{2} \omega_i^{(x_\tau,\tau)} E_+^{(x_\tau,\tau)}
+ \frac{1}{2} \omega_i^{(x_\tau,\tau)} E_-^{(x_\tau,\tau)}
\end{align*}
where $\omega^{(x_\tau,\tau)} := (x - x_\tau) / |x - x_\tau|$ is the analogue of $\omega$ for the origin $(x_\tau,\tau)$.  On the disk $D$ we have $\omega^{(x_\tau,\tau)} = -e_1 + O(\Lambda/R)$, and thus we informally have
$$ \Sigma ``\lesssim'' \dabb^{(x_\tau,\tau)} + \frac{\Lambda}{R} \nabla_{x,t}.$$
By many applications of the Leibnitz rule one thus obtains the rigorous estimates 
$$ |\nabla_{x,t} \Sigma^k \psi(0)| \lesssim \sum_{i+j \leq k} (\Lambda/R)^{k-i} |\nabla_{x,t}^{j+1} \tilde X^i \psi(0)|$$
on the disk $D$ with $\tilde X = \dabb^{(x_\tau,\tau)}$.  To obtain the regularity estimate \eqref{reg}, it thus suffices to show that
$$ \sum_{i+j \leq k} (\Lambda/R)^{k-i} 
E(\nabla_{x,t}^j \tilde X^i \psi[0]) \lesssim R^{100n} (\frac{\Lambda \mu}{R})^k.$$
Since $\psi$ has frequency $\mu$, it suffices by the reproducing formula \eqref{reproducing} to show that
$$ E(\tilde X^i \psi[0]) \lesssim R^{100n} (\frac{\Lambda \mu}{R})^i.$$
Observe that the ``good'' derivatives $\tilde X \lesssim \frac 1R \Gamma^{(x_\tau,\tau)}$.
Thus 
$$\aligned
E(\tilde X^i \psi[0]) &\lesssim (\frac \Lambda {R})^i \sum_{m\le i} 
\Lambda^{-m}E((\Gamma^{(x_\tau,\tau)})^m \psi[0])\\ & = 
 (\frac \Lambda {R})^i \sum_{m\le i} 
\Lambda^{-m}E((\Gamma^{(x_\tau,\tau)})^m \psi[\tau])\\ &= 
 (\frac \Lambda {R})^i \sum_{m\le i} 
\Lambda^{-m}E((\Gamma^{(x_\tau,\tau)})^m (\chi_{x_\tau} \phi[\tau]))
\endaligned
$$
The claim now follows from the Leibnitz rule, energy estimates for $\phi$, and
the frequency $\mu$ hypothesis.  This proves the localization estimate \eqref{local-2}.

Finally, we show the orthogonality estimate \eqref{ortho}.  By the Cotlar-Knapp-Stein lemma (see e.g. \cite{stein:large}) it suffices to show that
$$ |\langle \phi_T, \phi_{T'} \rangle_e| \lesssim (R + |x_0 - x'_0| + |x_\tau - x'_T|)^{-100n}$$
whenever $|x_0 - x'_0| + |x_\tau - x'_\tau| \gg r$ and $T = T(x_0,x_\tau)$, $T' = T(x'_0,x'_\tau)$ are in $\T$.

Fix $x_0, x'_0, x_\tau, x'_\tau$.  Define $x_t := x_0 + \frac{t}{\tau}(x_\tau-x_0)$ and $x'_t := x'_0 + \frac{t}{\tau}(x'_\tau - x'_0)$ for each $0 \leq t \leq \tau$.  Clearly we may find a time $R \leq t \leq 2R$ such that $|x_t - x'_t| \gg r$, which means that on this time slice the tubes $T$ and $T'$ are separated by $\gg r$.
Since $\langle \phi_T, \phi_{T'} \rangle_e = \langle \phi_T[t], \phi_{T'}[t] \rangle_e$, the claim follows from the localization estimate \eqref{local-2}.  This completes the second proof of Lemma \ref{tube-decomp}.

\section{Additional remarks}\label{remarks-sec}

$\bullet$ The induction on scales argument in Section \ref{induction-sec} lost a power of $R^\eps$ on a cube of length $R$.  This $R^\eps$ loss can be removed in a number of ways.  One possibility is to tighten the induction so that, if one inductively assumes a bilinear estimate with a constant of $A$ (instead of something like $R^\eps$), then one can deduce the same estimate with a constant of $A/2 + O(1)$.  By standard bootstrap arguments this gives the bound with $A = O(1)$.  This can be done provided one performs extremely careful tube decompositions in which one does not lose any constant at all (not even a factor of 2!) in the main term of the orthogonality inequality \eqref{ortho}, as well as the Bessel identity \eqref{bessel}.  This approach was carried out in \cite{tao:cone}, however that argument is heavily dependent on Plancherel's identity (among other things) and will probably not be extendable to rough metrics.

$\bullet$ Another option is to tighten the $L^2$ estimate to an $L^p$ estimate for some $p<2$ while still losing an $R^\eps$, in the spirit of \cite{wolff:cone}, \cite{tataru:bilinear}; this can be accomplished (for $p > (n+3)/(n+1)$) by using further geometric facts about the tubes\footnote{In our current argument, the only geometric fact used is that any two transverse tubes only intersect in a localized region of spacetime.  A further fact, used to good effect in \cite{wolff:cone} (see also \cite{tao:cone}) is that a tube which is transverse to an entire \emph{light cone} (not just another tube) will still only intersect in a localized region.  An alternative way to phrase this is that one can recover good energy estimates on light cones and similar null surfaces if one makes a transversality assumption.}. 
Standard ``epsilon-removal'' arguments (see \cite{borg:cone}, \cite{tv:cone2}; related arguments are in \cite{borg:kakeya}, \cite{tao:boch-rest}) exploiting some decay of the fundamental solution in a manner reminiscent of the Tomas-Stein or Strichartz arguments, can then be used to trade in the $R^\eps$ loss at the expense of $L^p$ integrability.

Since this argument does not appear to be well-known in the PDE literature we briefly sketch it as follows.  Suppose for instance that we have obtained an estimate of the form
$$
\| \phi \psi \|_{L^p(\cube_R)} \lesssim R^\eps E(\phi)^{1/2} E(\psi)^{1/2}
$$
for some $p<2$ and all $\eps > 0$, $R \gg 1$, and all free waves $\phi$, $\psi$ of frequency 1.  We sketch how this estimate implies the global estimate
\be{bil}
\| \phi \psi \|_{L^2_{x,t}} \lesssim E(\phi)^{1/2} E(\psi)^{1/2}.
\end{equation}
For simplicity we shall reduce to the case when $\phi_t = \pm i\sqrt{-\Delta} \phi$, and similarly for $\psi$.  In this case the energy inner product is equivalent to the usual inner product.  We shall now cease to distinguish between $\phi(t)$ and $\phi[t]$ to simplify the notation.

Firstly, we use finite speed of propagation to replace $\cube_R$ by a spacetime slab $\R^n \times I$, where $I$ is any interval of length $R$.  By Sobolev embedding we may then replace the $L^p_x$ norm with an $L^2_x$ norm:
\be{localized}
\| \phi \psi \|_{L^p_t L^2_x(\R^n \times I)} \lesssim R^\eps E(\phi)^{1/2} E(\psi)^{1/2}.
\end{equation}
Now we prove \eqref{bil}.  It suffices to prove the weak-type estimate
$$
\| \phi \psi \|_{L^{q,\infty}_t L^2_x} \lesssim E(\phi)^{1/2} E(\psi)^{1/2}
$$
for $q<p$.  By duality it suffices to show that
\be{dual-targ}
\int_\Omega \| \phi(t) \psi(t) \|_{L^2_x}\ dt \lesssim |\Omega|^{1/q}
E(\phi)^{1/2} E(\psi)^{1/2}.
\end{equation}
This estimate is clear for $|\Omega| \lesssim 1$ just by crudely bounding $\phi$ in $L^2$ and $\psi$ in $L^\infty$ (for instance), so we may assume $|\Omega| \gg 1$.

To compare with \eqref{dual-targ}, observe that the estimate \eqref{localized} yields
\be{near-miss}
\int_{\Omega \cap I} \| \phi(t) \psi(t) \|_{L^2_x}\ dt \lesssim |\Omega|^{1/p} |I|^\eps
E(\phi)^{1/2} E(\psi)^{1/2}
\end{equation}
for any interval $I$.  This is superior if $\Omega$ is concentrated in an interval $I$ of size $O(|\Omega|^C)$, but we still must consider the case when $\Omega$ is widely dispersed along the time axis.  To resolve this we use the decay of the fundamental solution as well as a $TT^*$ argument.

To see this, we fix $\Omega$ and let $R := |\Omega|^C$ for some large $C$ to be chosen later.  We then choose $\eps$ so that $|\Omega|^{1/q} \lesssim |\Omega|^{1/p} R^\eps$.  We then partition the time axis $\R$ into disjoint intervals $I$ of length $R$, and observe from \eqref{near-miss}, the hypothesis $q < p$, and Cauchy-Schwarz that we have the ``weakly global'' inequality
\be{weak-global}
\sum_I \int_{\Omega \cap I} \| \phi_I(t) \psi_I(t) \|_{L^2_x}\ dt \lesssim |\Omega|^{1/q} 
(\sum_I E(\phi_I))^{1/2} (\sum_I E(\psi_I))^{1/2}
\end{equation}
for all free waves $\phi_I$, $\psi_I$.  Our next objective is to bootstrap this inequality to a ``semi-weakly global'' inequality
\be{semiweak-global}
\sum_I \int_{\Omega \cap I} \| \phi(t) \psi_I(t) \|_{L^2_x}\ dt \lesssim |\Omega|^{1/q} 
E(\phi)^{1/2} (\sum_I E(\psi_I))^{1/2}.
\end{equation}
A similar bootstrap will then obtain \eqref{dual-targ}.

To obtain \eqref{semiweak-global}, we fix the $\psi_I$ and write $\phi[t] = U(t) P_1 \phi[0]$ (using the reproducing formula \eqref{reproducing}) to rewrite the estimate as 
$$
|\sum_I \int_{\Omega \cap I} \int U(t)P_1\phi(0) \psi_I(t) F(x,t) \ dt| \lesssim |\Omega|^{1/q} 
E(\phi)^{1/2} (\sum_I E(\psi_I))^{1/2} \| F \|_{L^\infty_t L^2_x}$$
By duality it suffices to show that
$$
E(\sum_I \int_{\Omega \cap I} P_1 U(-t) (\psi_I(t) F(t)) \ dt)^{1/2} \lesssim |\Omega|^{1/q} 
(\sum_I E(\psi_I))^{1/2} \| F \|_{L^\infty_t L^2_x}.$$
We square this, and expand this as
\be{tt*}
|\sum_{I,J} \int_{\Omega \cap I} \int_{\Omega \cap J} \langle P_1 U(-t) (\psi_I(t) F(t)), P_1 U(-s) (\psi_I(s) F(s)) \rangle\ ds dt|
\lesssim
(|\Omega|^{1/q} 
(\sum_I E(\psi_I))^{1/2} \| F \|_{L^\infty_t L^2_x})^2.
\end{equation}
On the other hand, if we start with \eqref{weak-global} and apply similar considerations we obtain
$$
|\sum_{I} \int_{\Omega \cap I} \int_{\Omega \cap I} \langle P_1 U(-t) (\psi_I(t) F(t)), P_1 U(-s) (\psi_I(s) F(s)) \rangle\ ds dt|
\lesssim
(|\Omega|^{1/q} 
(\sum_I E(\psi_I))^{1/2} \| F \|_{L^\infty_t L^2_x})^2;$$
in other words, we already control the diagonal contribution $I=J$ of \eqref{tt*}.  By Cauchy-Schwarz we also control the near-diagonal contribution when $I$ and $J$ are adjacent.  It remains to consider the contribution when $I$ and $J$ are not adjacent.  In this case $|t-s| \gtrsim R$.  To exploit this we rewrite this contribution as
$$
|\sum_{I,J} \int_{\Omega \cap I} \int_{\Omega \cap J} \langle P^2_1 U(s-t) (\psi_I(t) F(t)), \psi_J(s) F(s) \rangle\ ds dt|
$$
The decay estimate shows that $P^2_1 U(s-t)$ maps $L^1$ to $L^\infty$ with a norm of $O(R^{-(n-1)/2})$, thus we can estimate this by
$$
R^{-(n-1)/2} \sum_{I,J} \int_{\Omega \cap I} \int_{\Omega \cap J}  E(\psi_I) \| F\|_{L^\infty_t L^2_x} E(\psi_J) \| F \|_{L^\infty_t L^2_x}\ ds dt
$$
which we can easily bound by
$$ R^{-(n-1)/2} |\Omega|^2 ((\sum_I E(\psi_I))^{1/2} \| F \|_{L^\infty_t L^2_x})^2,$$
which is acceptable if we choose $R = |\Omega|^C$ sufficiently large.   (Note that we do not need the optimal decay of $(n-1)/2$ to run this argument).  This proves \eqref{semiweak-global}; the extension to \eqref{dual-targ} is done
 similarly and is omitted.  This completes the derivation of \eqref{bil}.

\medskip

$\bullet$ In what follows we sketch  a slightly different proof of 
Theorem \ref{null-1}, in which the parallel case is treated
 somewhat differently. We start, as before, with the tube
 decomposition
$$ \phi \approx \sum_{T \in \T} \phi_T ; 
\quad \psi \approx \sum_{T' \in \T'} \psi_{T'}. $$
We observe that the wave  packets $\phi_T$, $\psi_T$
verify   the properties \eqref{local} and 
\eqref{local-2} of definition \ref{packet-def}, not only
for $R/2\le t\le 5R/2$ but also at $t=0$.
This can easily be checked from either of the 
two construction of the wave packets. For simplicity
we shall assume that in fact, at $t=0$, both  $\phi_T,\psi_T$ are
localized to   balls of radius $r=R^{\f12+\eps}$.
As before we write:
$$
Q(\phi,\psi)\approx \sum_{T \in \T} \sum_{T' \in \T'} Q(\phi_T, \psi_{T'})
\approx Q_{||}+Q_{\times}
$$
with the transversal interaction component
$$
Q_{\times}=\sum_{\kappa, \kappa' \in \K: \kappa \sim \kappa'} Q( \phi_\kappa, \psi_{\kappa'} )
$$
and the parallel interaction component
$$
Q_{||}=\sum_{\kappa,\kappa' \in \K: 
\kappa \approx \kappa', l(\kappa) = l(\kappa') = \rho_0}
 Q( \phi_\kappa, \psi_{\kappa'} ).
$$
Recall that,
$$ \phi_\kappa := \sum_{T \in \T(\kappa)} \phi_T \hbox{ and }
\psi_\kappa := \sum_{T' \in \T'(\kappa)} \psi_{T'}.$$
Due to the additional localization  properties of the 
wave packets at $t=0$ we observe that $\phi_k(0),\psi_k(0)$
are  essentially localized in the cube $k$. We start by estimating
$Q_{||}$. Recalling, the estimate
$$ \| Q(\phi_k,\psi_{k'})(t) \|_2 \lesssim \frac{1}{R} 
(\| \Gamma \phi_k(t)\|_\infty \| \nabla_{x,t} \psi_{k'}(t)\|_2 +
 \|\nabla_{x,t} \phi_k\|_\infty \|\Gamma \psi_{k'}(t)\|_2),$$
we write
\beaa
\|Q_{||}\|_{L^2(\cube^{(0)}_R)}&\lesssim&\sum_{\kappa,\kappa' \in \K: 
\kappa \approx \kappa', l(\kappa) = l(\kappa') = \rho_0}
\| Q( \phi_\kappa, \psi_{\kappa'} )\|_{L^2(\cube^{(0)}_R)}\\
&\lesssim& R^{\f12}(\frac{\rho_0}{R})^{\frac{n+1}{2}}\sum_{\kappa,\kappa' \in \K: 
\kappa \approx \kappa', l(\kappa) = l(\kappa') = \rho_0}E(\phi_k(0))^{\f12}
 E(\psi_{k'}(0))^{\f12}\\
&\lesssim& R^{\f12}(\frac{\rho_0}{R})^{\frac{n+1}{2}}E(\phi)^{\f12}E(\psi)^{\f12}
\eeaa
We can therefore chose $\rho_0=R^{\frac{n}{n+1}}$ rather than 
$\rho_0\approx R^{\f12 +\eps}$ as used earlier.
The proof of the estimate for  $Q_{\times}$
proceeds exactly as before.

\section{Proof of Theorem \ref{null-2}}\label{null-sec}

Fix $\phi$, $\psi$, $\mu$, $\lambda$, $q$, $r$.  By a scaling argument we may take $\lambda = 1$.  We may assume that $\mu \lesssim 1$ since the left-hand side vanishes otherwise. Our task is now to prove 
\be{job-2}
\| P_\mu( \phi \psi ) \|_{L^{q/2}_t L^{r/2}_x(\R^{n+1})} \lesssim 
\mu^{ n - \frac{4}{q} - \frac{2n}{r} }
E(\phi)^{1/2} E(\psi)^{1/2}.
\end{equation}

If one simply discards the $P_\mu$ projection and uses H\"older and Strichartz then one does not obtain the decay in $\mu$ in \eqref{job-2}.  To recover this decay we will use the averaging property of $P_\mu$ at scale $1/\mu$ (given by \eqref{pointwise}), interacted against the fact that the fundamental solution at frequency 1 is essentially supported on a set of ``thickness'' 1.  (In this we are very much in the spirit of \cite{tao:lowreg}).

We divide space $\R^n$ into cubes $Q$ of side-length $1/\mu$.  From \eqref{p-split} we have
$$ \| P_\mu( \phi(t) \psi(t) ) \|_{ L^{r/2}_x(\R^n)} 
= (\sum_Q \| P_\mu( \phi(t) \psi(t) ) \|_{ L^{r/2}_x(Q)}^{r/2})^{2/r}.$$
By \eqref{pointwise} we have 
$$
\| P_\mu( \phi(t) \psi(t) ) \|_{ L^{r/2}_x(Q)}
\lesssim \mu^{n(1 - \frac{2}{r})} \sum_{Q'} (1 + \mu \dist(Q,Q'))^{-100n} 
\| \phi(t) \psi(t) \|_{ L^1_x(Q)}.$$
Inserting this into the previous and using Young's inequality, we obtain
$$ \| P_1( \phi(t) \psi(t) ) \|_{ L^{r/2}_x(\R^n)} 
\lesssim \mu^{n(1 - \frac{2}{r})} (\sum_Q \| \phi(t) \psi(t) \|_{ L^1_x(Q)}^{r/2})^{2/r}.$$
To prove \eqref{job-2}, it thus suffices by H\"older and Cauchy-Schwartz to show the ``improved Strichartz'' estimate
\be{improv-strichartz}
 \| (\sum_Q \| \phi(t) \|_{ L^2_x(Q)}^r)^{1/r} \|_{L^q_t}
\lesssim \mu^{-\frac{2}{q}} E(\phi)^{1/2}.
\end{equation}
together with the identical estimate for $\psi$.  

It suffices to prove \eqref{improv-strichartz} for $\phi$.  We remark that this estimate has essentially appeared in \cite{tao:lowreg}, and we shall repeat the proof here.  

Observe that the claim is trivial when $q = \infty$ (thanks to the embedding $l^r \subset l^2$), so we will assume that $q < \infty$.  This forces $r$ to be strictly greater than 2. 

We shall use the abstract Strichartz estimate in \cite{tao:keel} (Theorem 10.1).  Let $H$ be the Hilbert space of initial data $\phi[0]$ with frequency 1 with the energy inner product.  We let $U(t): H \to H$ be the map $U(t) \phi[0] := \phi[t]$; note that this map is unitary on $H$.  For any $0 \leq \theta \leq 1$, let $B_\theta$ denote the Banach space given by the norm
$$ \| \phi[t] \|_{B_\theta} := \mu^{-\theta \sigma} (\sum_Q \| \phi[t] \|_{ L^2_x(Q)}^{r_\theta})^{1/r_\theta}$$
where $r_\theta := 2/(1+\theta)$ and $\sigma := \frac{2r}{q(r-2)}$.  Since $(q,r)$ is admissible we observe that $0 < \sigma \leq (n-1)/2$.   Observe that the dual space $B^*_\theta$ is given by
$$ \| \phi[t] \|_{B^*_\theta} := \mu^{\theta \sigma} (\sum_Q \| \phi[t] \|_{ L^2_x(Q)}^{r_\theta'})^{1/r_\theta'}.$$
Also, one has the real interpolation embedding $(B^*_0, B^*_1)_{\theta,2} \subseteq B^*_\theta$ (see e.g. \cite{bergh:interp}).

To prove \eqref{improv-strichartz} it will suffice to prove the estimate
$$ \| U(t) \phi[0] \|_{L^q_t B^*_\theta} \lesssim \| \phi[0] \|_H$$
with $\theta := 1 - 2/r$ (so that $r_\theta' = r$).

Observe that $q = 2/\sigma\theta$.  By the abstract Strichartz estimate in \cite{tao:keel}, Theorem 10.1, we see that it will suffice to prove the decay estimate
$$ \| U(t) U(s)^* \phi[s] \|_{B^*_1} \lesssim |t-s|^{-\sigma}  \| \phi[s] \|_{B_1}$$
for all $s, t$.  From the group properties of the evolution we see that $U(t) U(s)^* = U(t-s)$, so by time translation invariance it suffices to show that
$$ \| U(t) \phi[0] \|_{B^*_1} \lesssim |t|^{-\sigma}  \| \phi[0] \|_{B_1}.$$
We expand out the norms to write this as
$$ \mu^{\sigma} \sup_Q \| \phi[t] \|_{L^2_x(Q)}
\lesssim \mu^{-\sigma} |t|^{-\sigma} \sum_Q \| \phi[0] \|_{L^2_x(Q)}.$$
It thus suffices to show that
$$ \| \phi[t] \|_{L^2_x(Q)}
\lesssim (1/\mu^2|t|)^\sigma \sum_{Q'} \| \phi[0] \|_{L^2_x(Q')}$$
for all cubes $Q$ of side-length $\mu$.

We may assume that $|t| \gg \mu^{-2}$ since the claim follows immediately from energy estimates otherwise.

We shall assume zero initial velocity $\phi_t(0) = 0$; the case of zero initial position is similar and is left to the reader.  
By the reproducing formula \eqref{reproducing} we may write $\phi(0) = \phi(0) * \psi(0)$ for some Schwartz function $\psi(0)$ with frequency 1.  We extend $\psi$ to be a free wave by setting zero initial velocity $\psi(0) = 0$.  By uniqueness we see that $\phi[t] = \phi(0) * \psi[t]$, thus $\psi[t]$ is a frequency-localized fundamental solution for the wave equation.  By Young's inequality we thus see that it suffices to show the decay estimate
$$  \| \psi[t] \|_{L^1_x(Q)}
\lesssim (1/\mu^2|t|)^\sigma$$
for all cubes $Q$ of side-length $1/\mu$.

Since $\psi[0]$ is a Schwartz function, we have the standard decay estimate 
$$ |\nabla^k_{x,t} \psi(x,t)| \lesssim |t|^{-(n-1)/2} (1 + ||t|-|x||)^{-k}$$
(see e.g. \cite{sogge:wave}; one can also use the machinery of Section \ref{vector-sec} to achieve this bound).  Since $\psi$ has frequency 1, we thus have in particular that
$$ |\psi[t](x)| \lesssim |t|^{-(n-1)/2} (1 + ||t|-|x||)^{-k}$$
for any $k$. A computation then shows that
$$  \| \psi[t] \|_{L^1_x(Q)}
\lesssim (1/\mu^2 |t|)^{(n-1)/2}.$$
Since $\sigma \leq (n-1)/2$, the claim follows.  This completes the proof of Theorem \ref{null-2}.

\section{Appendix: Proof of Corollaries \ref{cor-1}, \ref{cor-2}}\nonumber

We now prove Corollary \ref{cor-1}.  We first begin with the preliminary remark that in Theorem \ref{null-1} the spacetime cube $\cube_R$ can be replaced with the slab $[0,R] \times \R^n$.  To see this we divide into cases.

\textbf{Case 1: $R \min(\lambda,\mu) \geq 1$.}

In this case we use finite speed of propagation and smooth cutoff functions applied to the initial data to replace the energy $E(\phi)$ in \eqref{null-eq} by the more localized quantity
$$ \int (|\phi_t(x,0)|^2 + |\nabla \phi(x,0)|^2) (1 + \frac{\dist((x,0),\cube_R)}{R})^{-100n}\ dx.$$
(One has to construct the cutoff functions so as not to destroy the property of having frequency $\lambda$, but since $R \lambda \geq 1$ this can be done without violating the uncertainty principle.  See \cite{tao:cone} for further discussion).  Similarly for $E(\psi)$.  The claim then follows by partitioning $[0,R] \times \R^n$ into cubes $\cube_R$, applying the above localized estimates, and then Cauchy-Schwarz.

\textbf{Case 2: $R \min(\lambda,\mu) \leq 1$.}

Without loss of generality we may assume that $\lambda \leq \mu$.  From the definition of energy and Bernstein's inequality  we see that
$$ \| \nabla_{x,t} \phi \|_{L^\infty_t L^\infty_x} \lesssim \lambda^{n/2} E(\phi)^{1/2}$$
and
$$ \| \nabla_{x,t} \psi \|_{L^\infty_t L^2_x} \lesssim E(\psi)^{1/2}.$$
The claim then follows from  H\"older.

We now return to the proof of Corollary \ref{cor-1}.
Fix $Q$, $\phi$, $\psi$.  We use suitably chosen Littlewood-Paley  projections to split
$$
Q(\phi,\psi) = \sum_{\lambda, \mu} Q(P_\lambda \phi, P_\mu \psi).$$
We split into three cases: the high-low case $\lambda \gg \mu$, the high-high case $\lambda \sim \mu$, and the low-high case $\lambda \ll \mu$.

In the high-low case we see that $Q(P_\lambda \phi, P_\mu \psi)$ has Fourier transform supported in the region $|\xi| \sim \lambda$.  By orthogonality we thus have
$$ \| \sum_{\lambda \gg \mu} Q(P_\lambda \phi, P_\mu \psi) \|_{L^2_{t,x}([0,1] \times \R^n)}
\lesssim (\sum_{\lambda} 
\| \sum_{\mu: \mu \ll \lambda} Q(P_\lambda \phi, P_\mu \psi) \|_{L^2_{t,x}([0,1] \times \R^n)}^2)^{1/2}.$$
First consider the contribution when $\mu \lesssim 1$.  Applying Theorem \ref{null-1} (with the previously mentioned modification) we may estimate this by
$$
\lesssim 
(\sum_{\lambda} 
(\sum_{\mu: \mu \lesssim 1} \mu^{n/2}
E(P_\lambda \phi)^{1/2} E(P_\mu \psi)^{1/2})^2)^{1/2}.$$
This easily sums to be bounded by $E(\phi)^{1/2} E(\psi)^{1/2}$, which is acceptable.

Now consider the contribution when $\mu \gtrsim 1$.  Again using Theorem \ref{null-1}, we may bound this by
$$
\lesssim 
(\sum_{\lambda} 
(\sum_{\mu: 1 \lesssim \mu \ll \lambda} \mu^{(n-1)/2 + \eps}
E(P_\lambda \phi)^{1/2} E(P_\mu \psi)^{1/2})^2)^{1/2}.$$
This can be controlled by
$$
\lesssim 
(\sum_{\lambda} 
(\lambda^{(n-1)/2 + \eps}
E(P_\lambda \phi)^{1/2} E(\psi)^{1/2})^2)^{1/2}
\lesssim E( \langle D \rangle^{(n-1)/2+\eps} \phi)^{1/2} E(\psi)^{1/2}$$
as desired.

In the high-high case we use the triangle inequality followed by Theorem \ref{null-1}.  In the $\mu \lesssim 1$ case we thus obtain a bound of
$$ 
\sum_{\lambda \sim \mu \lesssim 1}
\mu^{n/2} E(P_\lambda \phi)^{1/2} E(P_\mu \psi)^{1/2} \lesssim E(\phi)^{1/2} E(\psi)^{1/2}$$
which is acceptable.  In the $\mu \gtrsim 1$ case we obtain a bound of
$$ 
\sum_{\lambda \sim \mu \gtrsim 1}
\lambda^{(n-1)/2 + \eps} E(P_\lambda \phi)^{1/2} E(P_\mu \psi)^{1/2}$$
which by Cauchy-Schwarz is bounded by $E(\langle D \rangle^{(n-1)/2+\eps} \phi)^{1/2} E(\psi)^{1/2}$ as desired.

The low-high case is similar to the high-low case.  We only show the treatment of the $\lambda \gtrsim 1$ contribution.  By orthogonality and Theorem \ref{null-1}, this is bounded by
$$
\lesssim 
(\sum_{\mu} 
(\sum_{\lambda: 1 \lesssim \lambda \ll \mu} \lambda^{(n-1)/2 + \eps}
E(P_\lambda \phi)^{1/2} E(P_\mu \psi)^{1/2})^2)^{1/2}
\lesssim (\sum_{\mu} 
E(\langle D \rangle^{(n-1)/2+2\eps} \phi)
E(P_\mu \psi))^{1/2}$$
which is acceptable (replacing $\eps$ by $\eps/2$ as necessary).  This concludes the proof of Corollary \ref{cor-1}.

We now prove Corollary \ref{cor-2}.  We use suitably chosen Littlewood-Paley  projections and the triangle inequality to split
$$ \||D|^{-\sigma}(\phi \psi) \|_{L^{q/2}_t L^{r/2}_x} \leq
 \sum_{\lambda, \mu} \| |D|^{-\sigma} (P_\lambda \phi P_\mu \psi) \|_{L^{q/2}_t L^{r/2}_x} $$
and consider the high-low, high-high, and low-low cases separately.

In the high-low case $\lambda \gg \mu$ the symbol $|D|^{-\sigma}$ can be replaced by $\lambda^{-\sigma}$.  We then use H\"older and Strichartz' inequality (see e.g. \cite{tao:keel}) to estimate
\begin{align*}
\| |D|^{-\sigma} (P_\lambda \phi P_\mu \psi) \|_{L^{q/2}_t L^{r/2}_x} 
&\lesssim \lambda^{-\sigma} \| P_\lambda \phi \|_{L^q_t L^r_x}
\| P_\mu \psi \|_{L^q_t L^r_x}\\
&\lesssim \lambda^{-\sigma}
\lambda^{\frac{n}{2} - \frac{1}{q} - \frac{n}{r}} E(P_\lambda \phi)^{1/2}
\mu^{\frac{n}{2} - \frac{1}{q} - \frac{n}{r}} E(P_\mu \psi)^{1/2}.
\end{align*}
Thus we may estimate this contribution by
$$ \sum_\lambda \lambda^{-\sigma}
\lambda^{\frac{n}{2} - \frac{1}{q} - \frac{n}{r}} E(P_\lambda \phi)^{1/2}
(\sum_{\mu \ll \lambda} \mu^{\frac{n}{2} - \frac{1}{q} - \frac{n}{r}} E(P_\mu \psi)^{1/2}).$$
The inner sum can be estimated by
$$
\sum_{\mu \ll \lambda} \mu^{\frac{n}{2} - \frac{1}{q} - \frac{n}{r}} E(P_\mu \psi)^{1/2}
\lesssim
\lambda^{\frac{\sigma}{2}}
E(|D|^{\frac{n-\sigma}{2} - \frac{1}{q} - \frac{n}{r}} \psi)^{1/2},$$
and the claim then follows by performing the $\lambda$ summation.

The low-high case $\lambda \ll \mu$ is similar, so we turn to the high-high case $\lambda \sim \mu$.  To simplify the notation we shall take $\lambda = \mu$, although this does not substantially affect the argument.  By Cauchy-Schwarz it suffices to show that
$$
\| |D|^{-\sigma} (P_\lambda \phi P_\lambda \psi) \|_{L^{q/2}_t L^{r/2}_x} 
\lesssim 
E(|D|^{\frac{n-\sigma}{2} - \frac{1}{q} - \frac{n}{r}} P_\lambda \phi)^{1/2}
E(|D|^{\frac{n-\sigma}{2} - \frac{1}{q} - \frac{n}{r}} P_\lambda \psi)^{1/2}.$$
From scale invariance we may take $\lambda \sim 1$, so that the derivatives $|D|$ on the right hand side may be discarded.  

The function $P_\lambda \phi P_\lambda \psi$ has frequency $\lesssim 1$.  By the Littlewood-Paley  decomposition and the triangle inequality it thus suffices to show that
$$ 
\| P_\nu |D|^{-\sigma} ( P_\lambda \phi P_\lambda \psi ) \|_{L^{q/2}_t L^{r/2}_x} \lesssim \nu^\eps E(P_\lambda \phi)^{1/2} E(P_\lambda \psi)^{1/2}$$
for some $\eps > 0$.  But we may replace $|D|^{-\sigma}$ by $\nu^{-\sigma}$, and the claim then follows from  Theorem \ref{null-2} and the hypothesis $\sigma < n - \frac{4}{q} - \frac{2n}{r}$.

\end{document}